\documentclass[11pt]{article}
\usepackage[dvips]{graphics}
\input{amssym.def}
\input{amssym.tex}

\title{\Large\bf SINGULAR INTEGRAL EQUATIONS\\ WITH TWO FIXED SINGULARITIES
AND APPLICATIONS TO FRACTURED COMPOSITES}
\author{Y.A. ANTIPOV\\
antipov@math.lsu.edu\\ 
Department of Mathematics, Louisiana State University\\
Baton Rouge LA 70803, USA}

\newcommand{\bfm}[1]{\mbox{\boldmath ${#1}$}}

\newcommand{\beqa}{\begin{eqnarray}}
\newcommand{\eeqa}[1]{\label{#1}\end{eqnarray}}
\newcommand{\bequ}{\begin{equation}}
\newcommand{\eequ}[1]{\label{#1}\end{equation}}

\newcommand{\Md}{\partial}

\newcommand{\ov}[1]{\overline{#1}}

\newcommand{\Ga}{\alpha}
\newcommand{\Gb}{\beta}
\newcommand{\Gd}{\delta}

\newcommand{\Gve}{\varepsilon}          
\newcommand{\Gf}{\phi}
\newcommand{\Gvf}{\varphi}
\newcommand{\Gg}{\gamma}
\newcommand{\Gc}{\chi}

\newcommand{\Gl}{\lambda}
\newcommand{\Gn}{\eta}

\newcommand{\Gr}{\rho}

\newcommand{\Gs}{\sigma}

\newcommand{\Go}{\omega}
\newcommand{\Gx}{\xi}
\newcommand{\Gy}{\psi}
\newcommand{\Gz}{\zeta}
\newcommand{\GD}{\Delta}
\newcommand{\GF}{\Phi}
\newcommand{\GG}{\Gamma}
\newcommand{\GL}{\Lambda}
\newcommand{\GP}{\Pi}

\newcommand{\GO}{\Omega}

\newcommand{\BGF}{\bfm\Phi}

\newcommand{\BGY}{\bfm\Psi}

\newcommand{\CK}{{\cal K}}

\newcommand{\CS}{{\cal S}}

\def\Bb{{\bf b}}

\def\Bg{{\bf g}}

\def\BJ{{\bf J}}

\def\BV{{\bf V}}

\newcommand{\beq}{\begin{equation}}
\newcommand{\eeq}{\end{equation}}
\newcommand{\barr}{\begin{eqnarray}}
\newcommand{\earr}{\end{eqnarray}}
\newcommand{\beqn}{\begin{equation*}}
\newcommand{\eeqn}{\end{equation*}}
\newcommand{\barrn}{\begin{eqnarray*}}
\newcommand{\earrn}{\end{eqnarray*}}

\newcommand{\fr}{\frac}

\newcommand{\diag}{\mbox{diag}}

\newcommand{\supp}{\mathop{\rm supp}\nolimits}

\newcommand{\sgn}{\mathop{\rm sgn}\nolimits}

\newcommand{\I}{\mathop{\rm Im}\nolimits}
\newcommand{\R}{\mathop{\rm Re}\nolimits}
\newcommand{\const}{\mbox{const}}

\begin{document}
\maketitle

\begin{abstract}

A symmetric characteristic singular integral equation with two fixed singularities at the endpoints in the class
of  functions bounded at the ends  is analyzed.  It reduces to a vector Hilbert problem for a half-disc
and then to a vector Riemann-Hilbert problem on a real axis with a piecewise constant matrix coefficient
that has two points of discontinuity.
A condition of solvability and a closed-form solution to the integral equation are derived. 
For the Chebyshev polynomials of the first kind in the right hand-side, the solution of the
integral equation is expressed in terms of two nonorthogonal polynomials with 
associated weights. 
Based on this new generalized spectral relation for the singular operator with two fixed singularities 
an approximate solution 
to the complete singular integral equation is derived by
recasting it as an infinite system of linear algebraic equations of the second kind.
The method is illustrated by solving two problems of fracture mechanics,
the antiplane and plane strain problems for a finite crack in a composite plane. The plane is formed
by a strip and two half-planes; the elastic constants of the strip are different from those of the half-planes.
The crack is orthogonal to the interfaces, and it is located in the strip with the ends lying
in the interfaces. Numerical results are reported and discussed.
 
\end{abstract}

\setcounter{equation}{0}

\section{Introduction}

The method of orthogonal polynomials, a particular realization of the general Bubnov-Galerkin
method, has been successfully applied for singular integral equations since publication
of the work by Klubin  {\bf(\ref{klu})} who employed the spectral properties of the logarithmic kernel $\ln|x-\Gx|$
and the Weber-Schafheitlin integral $W(x,\Gx)=\int_0^\infty J_0(tx)J_0(t\xi)dt$, 
$J_0(y)$ is the Bessel function, to obtain series representations
of the solutions to the corresponding singular integral equations. In particular, by making use of the
spectral relation 
\beq
\int_0^1 W(x,\Gx)\fr{\Gx P_{2n}(\sqrt{1-\Gx^2})d\Gx}{\sqrt{1-\Gx^2}}=\fr{\pi}{2}
\left[\fr{(2n-1)!!}{(2n)!!]}\right]^2P_{2n}(\sqrt{1-x^2}), \quad 0<x<1,\quad n=0,1,\ldots,
\label{0.1}
\eeq
where $P_m(x)$ are the Legendre polynomials, Klubin derived an efficient solution to the contact problem on a circular plate lying on an elastic foundation.
The method of orthogonal polynomials was further developed and employed by many researchers including  {\bf(\ref{pop1})},  {\bf(\ref{pop2})},  {\bf(\ref{erd})}.
This scheme applied to a singular integral equation $\int_a^b M(x,\Gx)\Gvf(\Gx)d\Gx=f(x)$, $a<x<b$, requires to determine the singularities of the solution at the endpoints in the class prescribed and represent the kernel as $M(x,\Gx)=\GP(x,\Gx)+K(x,\Gx)$.
The function $\GP(x,\Gx)$ is the dominant singular kernel, while the second part
is normally bounded or may have a weaker singularity as $x=\Gx$. The method can be
successfully applied if $\GP(x,\Gx)$ is a polynomial kernel   {\bf(\ref{pop2})}, that is a function
 satisfying the spectral relations
\beq
\int_a^b \left(\begin{array}{c}\GP(x,\Gx)\\
\GP(\Gx,x)\\
\end{array}
\right) p_\pm(\Gx)\pi_n^\pm(\Gx)d\Gx=\Gs_n g_\pm(x)\pi_n^\mp(x), \quad a<x<b,\quad n=0,1,\ldots,
\label{0.2}
\eeq
where  $\Gs_n\ne 0$ and $\pi_n^\pm(x)$ are orthonormal polynomials 
with the weights $w_\pm(x)=p_\pm(x) g_\mp(x)$ in the segment $(a,b)$.
Employing  the system of functions 
$p_+(x)\pi_n^+(x)$ ($n=1,2,\ldots$) as basis functions, expanding 
the unknown function as 
\beq
\Gvf(x)=p_+(x)\sum_{m=0}^\infty a_m \pi_m^+(x)
\label{0.4}
\eeq
and inserting the series into the integral equation enable us to reduce the equation to
\beq
\Gs_na_n+\sum_{m=0}^\infty d_{nm} a_m=f_n, \quad n=0,1,\ldots,
\label{0.5}
\eeq
an infinite system of linear algebraic equations of the second kind. Under certain 
conditions, normally satisfied in applications, it is possible to prove the convergence of an 
approximate solution to the exact one.
An approximate solution to the infinite system (\ref{0.5}) is derived by the reduction method.

Many problems of thin plate theory  {\bf(\ref{ant})} and contact and fracture mechanics  {\bf(\ref{moi})}
are governed
by singular singular integral equations of the form 
\beq
\int_0^1[S(x,\Gx)+K(x,\Gx)]\Gvf(\Gx)d\Gx=f(x), \quad 0<x<1,
\label{0.6}
\eeq
where $S(x,\Gx)$ is the sum of the Cauchy kernel and the kernel whose only singularities
are at the endpoints $x=\Gx=0$ and $x=\Gx=1$.  Motivated by the antiplane problem for a crack
in a strip placed between two half-planes of different shear moduli when the crack is orthogonal 
to the  interfaces, Moiseyev and Popov  {\bf(\ref{moi})} analyzed the case when
\beq
S(x,\Gx)=\fr{1}{\Gx-x}+\fr{\Gb_0}{\Gx+x}+\fr{\Gb_1}{\Gx+x-2}+R_0(x,\Gx), 
\label{0.7}
\eeq
 $\Gb_0$ and $\Gb_1$ are real,   $|\Gb_j|<1$, and $R_0(x,\Gx)$ is a specifically
chosen regular kernel. They reduced the problem to the vector Riemann-Hilbert problem with a 
piecewise constant matrix coefficient with three points
of discontinuity and solved it in terms of some quadratures and  the hypergeometric functions.
They did not find a spectral relation for the singular operator
with the fixed singularities. The  approximate scheme outlined was not tested and found to be hard 
to implement.

The main goal of this paper is to derive an exact solution  to equation (\ref{0.6})  in a simple form
when $K(x,\xi)\equiv 0$,   $2S(x,\Gx)=\cot\fr12\pi(\Gx-x)+\Gb\cot\fr12\pi(\Gx+x)$,
 $\Gb_0=\Gb_1=\Gb$, $\Gb\in(-\infty, +\infty)$, and develop a 
sufficient numerical scheme for the complete singular integral  equation with an arbitrary
kernel $K(x,\Gx)$ having at most weak singularities in the line $x=\Gx$ and weak fixed singularities
at the endpoints.
Specifically, we aim to derive a spectral  relation of the form
\beq
\int_0^1\left[\cot\fr{\pi(\Gx-x)}{2}+\Gb\cot\fr{\pi(\Gx+x)}{2}\right]\Gf_n(\Gx)d\Gx=
\Gs_n g(x)\pi_{n}(x), 
\quad
 0<x<1, \quad n=0,1,\ldots,
\label{0.8}
\eeq
and use it as the core of an approximate scheme for equation (\ref{0.6}).
Here,  $\Gf_n(x)=p^{(1)}(x)\pi^{(1)}_n(x)
+p^{(2)}(x)\pi^{(2)}_n(x)$,  
$\pi_n^{(j)}(x)$ ($j=1,2$) are some degree-$n$ trigonometric
polynomials not necessary orthogonal, $p^{(j)}(x)$ 
are their weights, and
$\pi_n(x)$ are 
degree-$n$ orthogonal trigonometric polynomials, $\Gs_n\ne 0$, $g(x)>0$.

The paper is organized as follows. In section 2, we analyze the characteristic singular integral
equation (\ref{2.1}) by employing the method of the vector Hilbert problem for a half-disc 
proposed in {\bf(\ref{moi})}. The characteristic equation (\ref{2.1})
considered in the present paper in the class of functions bounded at the endpoints leads to a vector Riemann-Hilbert problem
whose matrix coefficient has two points of discontinuity. This enables us to write down the solvability condition and the solution to the integral equation explicitly. Motivated by possible applications that may arise in future applications
we consider all possible cases for the major parameter 
$\Gb$ of the kernel, $\Gb\in(-\infty,\infty)$, not only when $|\Gb|<1$. 
In the case $\Gb=0$, the kernel becomes the Hilbert kernel, and we show that
our general formulas for the solvability condition and the solution reduce to the ones consistent with the known
results {\bf(\ref{gak})}, p.426.

In section 3, we consider the complete singular integral equation with two fixed singularities in the class of 
functions bounded at the ends when the right-hand side is defined up to an arbitrary constant.
First we derive the relation (\ref{0.8}) with the right-hand side 
chosen to be the Chebyshev trigonometric polynomial. 
We further employ these new spectral relations to derive an efficient approximate solution 
to the complete singular integral equation with two fixed singularities by converting
the integral equation into an infinite system of linear algebraic equations of the second kind.

In section 4, the method is tested by solving the antiplane problem on a finite crack 
orthogonal to the interfaces between a strip and two half-planes when
the shear moduli of the half-planes are the same, while the shear modulus of the strip
is different, and the crack lies in the strip. Section 5 generalizes the method for 
the singular integral equation with two fixed singularities (\ref{5.9}) that governs the 
corresponding
plane strain problem. In this case we analyze the singularities of the solution
at the endpoints and construct a new singular integral operator  associated 
with (\ref{5.9}) and whose spectral properties are studied in Section 3. A quick numerical test
applied to equation (\ref{5.18}) confirms the efficiency of the method
for plane problems as well.

  In appendix A, we adjust to our case
 the proof {\bf(\ref{moi})} of the equivalence of the singular integral equation 
and the vector Hilbert problem for two analytic functions in the half-disc. In appendix B we show
that the cases $\Gb=\pm 1$ reduce to an equation whose exact solution can be obtained by the Hilbert inversion formula
{\bf(\ref{gak})}, p.244; {\bf(\ref{mus})}, p.69.
Appendix C computes certain auxiliary integrals needed for the derivation of the spectral relation. Appendix D employs the method of classical orthogonal polynomials
for an approximate solution of the complete singular integral equation with the Cauchy kernel in the class of functions vanishing at the endpoints.

\setcounter{equation}{0}

\section{Characteristic singular integral equation with fixed singularities}\label{s2}

In this section, we aim to construct the exact solution 
to the singular integral equation
\beq
\fr12\int_0^1
\left[\cot\fr{\pi(\Gx-x)}{2}+\Gb\cot\fr{\pi(\Gx+x)}{2}
\right]\Gf(\Gx)d\Gx=f(x), \quad 0<x<1,
\label{2.1}
\eeq
in the class of functions bounded at the points $x=0$ and $x=1$ and
 H\"older-continuous in the interval $(0,1)$, $\Gf(x)\in H(0,1)$.  Here, $\Gb$ is a real parameter, and $f(x)\in H[0,1]$. Denote the kernel of the equation as
\beq 
S(x,\Gx)=\fr12\cot\fr{\pi(\Gx-x)}{2}+\fr{\Gb}{2}\cot\fr{\pi(\Gx+x)}{2}.
\label{2.2}
\eeq
The first term of the kernel has a singularity in the line $\Gx=x$, $x\in[0,1]$, while the second term has fixed singularities at the endpoints $\Gx=x=0$
and  $\Gx=x=1$. The kernel admits the representation
\beq
S(x,\Gx)=\fr{1}{\pi}\left(\fr{1}{\Gx-x}+\fr{\Gb}{\Gx+x}+\fr{\Gb}{\Gx+x-2}\right)+R_0(x,\Gx),
\label{2.3}
\eeq
where $R_0(x,\Gx)$ is a regular kernel. We call (\ref{2.1}) the characteristic equation to distinguish it from  the complete
singular integral equation to be analyzed in section \ref{s3}.

\subsection{Vector Hilbert and Riemann-Hilbert problems associated with the integral equation}\label{s2.1}

To construct an exact solution to the characteristic integral equation (\ref{2.1}), first we
transform it into an equation on the upper arc $L$ of the unit circle centered at the origin with the starting and terminal points 1 and -1, respectively. 
Let $t=e^{i\pi x}$, $\tau=e^{i\pi\Gx}$, $\Gf(x)=u(e^{i\pi x})$, and $f(x)=v(e^{i\pi x})$. Then (\ref{2.1}) becomes
\beq
\fr{1}{2\pi}\int_L\left[\fr{\tau+t}{\tau-t}-\fr{\Gb(1+\tau t)}{1-\tau t}\right]\fr{u(\tau)d\tau}{\tau}=v(t),\quad t\in L.
\label{2.4}
\eeq
According to the theorem to be stated below this integral equation  is equivalent to a certain vector Hilbert boundary value problem for a half-disc. 
\vspace{.1in}

Theorem 2.1.  Let
$\Gvf_1(z)$ and $\Gvf_2(z)$ be two functions  analytic in the 
upper half-disc $D=\{z\in{\Bbb C}: |z|<1, \I z>0\}$, 
H\"older-continuous up to the boundary $\Md D=L\cup(-1,1)$, bounded
at the points $z=\pm 1$ and  satisfying the 
boundary conditions
$$
\R[\Gvf_1(t)-\Gvf_2(t)]=0, \quad \I[\Gb_0\Gvf_1(t)+\Gvf_2(t)]=0,\quad -1<t<1,\quad \Gb_0=\fr{1+\Gb}{1-\Gb},
$$
\beq
\R\Gvf_1(t)=0,\quad \R\Gvf_2(t)=u(t), \quad t\in L,
\label{2.9}
\eeq
and the additional condition 
\beq
\lim_{z\to 0}[\Gvf_1(z)-\Gvf_2(z)]=0.
\label{2.10}
\eeq
Then the function $u(t)$, $t\in L$, solves the singular integral equation (\ref{2.4})
with $v(t)=-\I\Gvf_2(t)$
in the class of functions $H(L)$ bounded at the points $t=\pm 1$.

Conversely, 
let $u(t)$ be a solution to the integral equation (\ref{2.4}) in the class of functions H\"older-continuous everywhere in the contour $L$  and bounded at 
the ending points.  Denote
\beq
\Gg(z)=\fr{\Gb}{2\pi i}[\ln z-\ln(-z)], \quad -\pi<\arg z<\pi.
\label{2.7}
\eeq
Then the functions
$$
\Gvf_1(z)=\fr{1}{2\pi i}\int_L
\left\{[\Gg(z)-\Gg(\tau)]\fr{\tau+z}{\tau-z}+[1-\Gg(z)-\Gg(\tau)]\fr{1+\tau z}{1-\tau z}\right\}\fr{u(\tau)d\tau}{\tau},
$$
\beq
\Gvf_2(z)=\fr{1}{2\pi i}\int_L
\left\{[1+\Gg(z)-\Gg(\tau)]\fr{\tau+z}{\tau-z}-[\Gg(z)+\Gg(\tau)]\fr{1+\tau z}{1-\tau z}\right\}\fr{u(\tau)d\tau}{\tau}
\label{2.8}
\eeq
form the solution of the Hilbert boundary value problem 
$$
\R[\Gvf_1(t)-\Gvf_2(t)]=0, \quad \I[\Gb_0\Gvf_1(t)+\Gvf_2(t)]=0,\quad -1<t<1,
$$
\beq
\R\Gvf_1(t)=0,\quad \I\Gvf_2(t)=-v(t), \quad t\in L,
\label{2.9'}
\eeq 
in the class of functions bounded at the points $z=\pm 1$, satisfying the condition (\ref{2.10}), and
$\R\Gvf_2(t)=u(t)$, $t\in L$.

\vspace{.1in}

This theorem follows from the results derived  for a more general case in {\bf(\ref{moi})}.
A proof adjusted for the case under consideration is presented in
appendix A.

We now proceed to the construction of the solution to the Hilbert boundary value problem (\ref{2.9}), (\ref{2.10}). First we conformally map
the domain $D$ onto the lower half-plane $\I s(z)<0$,
\beq
s=\left(\fr{z-1}{z+1}\right)^2, \quad z=\fr{1+s^{1/2}}{1-s^{1/2}}, \quad z\in D,\quad  \pi<\arg s<2\pi.
\label{2.11}
\eeq
In this way the contour $L$ is mapped onto the negative semi-axis $-\infty<s<0$, while the segment $(-1,1)$
is mapped onto the positive semi-axis $0<s<+\infty$. The points $z=1$ and $z=-1$ fall onto the points $s=0$
and $s=\infty$, respectively. Next we define the vector
\beq
\BGF(s)=\left(\begin{array}{c}
\Gvf_1(z(s))\\
\Gvf_2(z(s))\\
\end{array}
\right),\quad \I s<0,
\label{2.12}
\eeq
analytic in the lower $s$-half-plane and extend its definition into the upper half-plane by the symmetry law
\beq
\BGF(s)=\diag\{-1,1\}\ov{\BGF(\ov{s})}, \quad \I s>0.
\label{2.13}
\eeq
On denoting
\beq
\BGF^-(\Gs)=\BGF(\Gs-i0), \quad \BGF^+(\Gs)=\diag\{-1,1\}
\ov{\BGF(\Gs-i0)},\quad -\infty<\Gs<+\infty,
\label{2.14} 
\eeq
we write the Hilbert problem (\ref{2.9}) in the form of the vector Riemann-Hilbert 
problem with a piece-wise constant matrix coefficient
\beq
\BGF^+(\Gs)=G(\Gs)\BGF^-(\Gs)+\Bg(\Gs), \quad -\infty<\Gs<+\infty,
\label{2.15}
\eeq
where
$$
G(\Gs)=\left\{
\begin{array}{cc}
I, & -\infty<\Gs<0,\\
G_0, & 0<\Gs<+\infty,\\
\end{array}
\right.
\quad
G_0=\left(
\begin{array}{cc}
-\Gb &  \Gb-1\\
\Gb+1 & -\Gb\\
\end{array}
\right),\quad I=\diag\{1,1\},
$$
\beq
\Bg(\Gs)=\left\{
\begin{array}{cc}
2iv(z(\Gs))\BJ, & -\infty<\Gs<0,\\
\bf{0}, & 0<\Gs<+\infty,\\
\end{array}
\right.
\quad
\BJ=\left(
\begin{array}{c}
0\\
1\\
\end{array}
\right).
\label{2.16}
\eeq

\subsection{Case $|\Gb|<1$}\label{s2.2}

In this case the matrix $G_0$ has complex-conjugate eigenvalues, 
$\Gl_1=-\Gb+i\sqrt{1-\Gb^2}$ and $\Gl_2=-\Gb-i\sqrt{1-\Gb^2}$, and admits
the splitting $G_0=T \GL_0 T^{-1}$, where
\beq
T=\left(
\begin{array}{cc}
1 &  1\\
-i\sqrt{\Gb_0} & i\sqrt{\Gb_0}\\
\end{array}
\right),\quad \GL_0=\diag\{\Gl_1,\Gl_2\},
\label{2.19}
 \eeq
and $\Gb_0$ is the parameter defined in (\ref{2.9}).
It is natural now to introduce the following vectors and matrices:
\beq
\BGF_0(s)=T^{-1}\BGF(s), \quad
\GL(\Gs)=\left\{
\begin{array}{cc}
I, & -\infty<\Gs<0,\\
\GL_0, & 0<\Gs<+\infty,\\
\end{array}
\right.\quad \Bg_0(\Gs)=T^{-1}\Bg(\Gs).
\label{2.20}
\eeq
The  vector Riemann-Hilbert problem has been decoupled, and the one-sided limits
$\BGF_0^\pm(\Gs)$ of the new vector $\BGF_0(s)$ satisfy the boundary condition
\beq
\BGF^+_0(\Gs)=\GL(\Gs)\BGF^-_0(\Gs)+\Bg_0(\Gs), \quad -\infty<\Gs<+\infty.
\label{2.21}
\eeq
The diagonal entries of the matrix $\GL(\Gs)$ can be factorized by means of the Cauchy integral
\beq
\Gl_j=\fr{\Gc_j^+(\Gs)}{\Gc^-_j(\Gs)}, \quad \Gs\in(0,+\infty),\quad j=1,2,
\label{2.22} 
\eeq
where $\Gc_j^\pm(\Gs)$ are the one-sided limits as $s\to\Gs\pm i0$ of the function
\beq
\Gc_j(s)=\exp\left\{\fr{1}{2\pi i}
\int_0^\infty
\left(\fr{1}{\Gs-s}-\fr{1}{\Gs-s_0}\right)\log\Gl_jd\Gs\right\}
=C_j^0s^{-(2\pi i)^{-1}\log\Gl_j}, 
\label{2.22'}
\eeq
and
\beq
 -2\pi\le\arg\Gl_j\le 0.
\label{2.22.0}
\eeq
The choice of the range for the argument of the eigenvalues is due to the class of solutions.
Here, $s_0\in {\Bbb C}\setminus(0,+\infty)$ is an arbitrary fixed point, and $C_j^0$ $(j=1,2)$ are constants.
 Denote
\beq
\Gr_j=-\fr{\log\Gl_j}{2\pi i}, \quad \Gd=\fr{1}{\pi}\tan^{-1}\fr{\sqrt{1-\Gb^2}}{\Gb}\in\left(-\fr{1}{2},\fr{1}{2}\right),
\quad j=1,2.
\label{2.22.1}
\eeq
Since $|\Gl_j|=1$, $j=1,2$, and because of the range for $\arg \Gl_j$ defined by (\ref{2.22.0}),
 we find
\beq
\Gr_1=\left\{
\begin{array}{cc}
\fr12+\fr{\Gd}{2}, & \Gb>0,\\
1+\fr{\Gd}{2}, & \Gb<0,\\
\end{array}
\right.
\quad 
\Gr_2=\left\{
\begin{array}{cc}
\fr12-\fr{\Gd}{2}, & \Gb>0,\\
-\fr{\Gd}{2}, & \Gb<0.\\
\end{array}
\right.
\label{2.22.2}
\eeq
We consequently derive the factorization of the diagonal matrix $\GL(\Gs)=X^+(\Gs)[X^-(\Gs)]^{-1}$,
$\Gs\in(-\infty,+\infty)$, where $X^\pm(\Gs)$ are the one-sided limits as $s\to \Gs\pm i0$ of the matrix $X(s)$
\beq
X(s)=\diag\{s^{\Gr_1},s^{\Gr_2}\}.
\label{2.22.3} 
\eeq
The branches of the functions $s^{\Gr_1}$ and $s^{\Gr_2}$ are fixed by cutting 
the $s$-plane along the positive semi-axis and selecting  $\arg  s\in[0,2\pi]$.
Next, by the reasoning usual in the theory of the Riemann-Hilbert problem we obtain
\beq
\BGF_0(s)=X(s)\BV(s), \quad \BV(s)=\fr{1}{2\pi i\sqrt{\Gb_0}}\int_{-\infty}^0\fr{v(z(\Gs))}{\Gs-s}\left(
\begin{array}{c}
-\Gs^{-\Gr_1}\\
\Gs^{-\Gr_2}\\
\end{array}\right) d\Gs.
\label{2.23}
\eeq	  
On returning to the $z$-plane by means of (\ref{2.11}), we transform the vector-function $\BV(s)$ as 
\beq
\BV(s)=\fr{(z+1)^2}{2\pi i\sqrt{\Gb_0}}\int_L\fr{v(t)}{(t-z)(1-tz)}\left(
\begin{array}{c}
-\Go_1(t)\\
\Go_2(t)\\
\end{array}\right) dt.
\label{2.24}
\eeq	  
Here,
\beq
\Go_j(z)=\left(\fr{z-1}{z+1}\right)^{1-2\Gr_j}, \quad j=1,2,
\label{2.25}
\eeq
and the branch of $\Go_j(z)$ in the $z$-plane cut along the line joining the points 
$z=1$ and $z=-1$ and passing through the infinite point
is chosen such that $\Go_j(0)=-e^{-2\pi i\Gr_j}$, $j=1,2$.
On employing formulas (\ref{2.20}) and (\ref{2.12})  it is now easy to derive from here representations for the
functions $\Gvf_1(z)$ and $\Gvf_2(z)$, the solution of the Hilbert problem (\ref{2.9}),
$$
\Gvf_1(z)=
\fr{z^2-1}{2\pi i\sqrt{\Gb_0}}\int_L
\left[-\fr{\Go_1(t)}{\Go_1(z)}+\fr{\Go_2(t)}{\Go_2(z)}\right]\fr{v(t)dt}{(t-z)(1-tz)},
$$
\beq
\Gvf_2(z)=\fr{z^2-1}{2\pi}\int_L\left[\fr{\Go_1(t)}{\Go_1(z)}+\fr{\Go_2(t)}{\Go_2(z)}\right]\fr{v(t)dt}{(t-z)(1-tz)}.
\label{2.26}
\eeq
These functions have to satisfy the necessary and sufficient condition (\ref{2.10}) for the equivalence
of the Hilbert problem (\ref{2.9})
and the integral equation (\ref{2.4}). It reads
\beq
\int_L\left[\left(\fr{1}{\sqrt{\Gb_0}}+i\right)\fr{\Go_1(t)}{\Go_1(0)}-\left(\fr{1}{\sqrt{\Gb_0}}-i\right)\fr{\Go_2(t)}{\Go_2(0)}\right
]\fr{v(t)dt}{t}=0.
\label{2.27}
\eeq
If this condition is satisfied, then the solution to equation (\ref{2.4}) is expressed
through the solution to the Hilbert problem by the formula $u(t)=\R\Gvf_2(t)$, $t\in L$.

In order to recover the solution to the original characteristic equation (\ref{2.1}) and verify the boundary conditions on the contour $L$, we make
the reverse substitution $t=e^{i\pi x}$, $\tau=e^{i\pi\Gx}$ and utilize the Sokhotski-Plemelj formulas. Since
\beq
\fr{(t^2-1)d\tau}{(\tau-t)(1-\tau t)}=\fr{\pi\sin\pi x d\Gx}{\cos\pi\Gx-\cos\pi x},
\label{2.28}
\eeq
we obtain
$$
\Gvf_1(t)=\fr{\sin\pi x}{2i\sqrt{\Gb_0}}\int_0^1\fr{v(e^{i\pi\Gx})[a^{2\Gr_1-1}(x,\Gx)-a^{-2\Gr_1+1}(x,\Gx)]d\Gx}{\cos\pi\Gx-\cos\pi x},
$$
\beq
\Gvf_2(t)=-iv(t)+\fr{\sin\pi x}{2}\int_0^1\fr{v(e^{i\pi\Gx})[a^{2\Gr_1-1}(x,\Gx)+a^{-2\Gr_1+1}(x,\Gx)]d\Gx}{\cos\pi\Gx-\cos\pi x}. \quad t\in L, \quad 0<x<1,
\label{2.29}
\eeq
where
\beq
a(x,\Gx)=\tan\fr{\pi \Gx}{2}\cot\fr{\pi x}{2}.
\label{2.30}
\eeq
These formulas imply
$$
\R\Gvf_1(t)=0, \quad \I\Gvf_2(t)=-v(t), \quad t\in L,
$$
\beq
\R\Gvf_2(t)=\fr{\sin\pi x}{2}\int_0^1\fr{v(e^{i\pi\Gx})[a^{2\Gr_1-1}(x,\Gx)+a^{-2\Gr_1+1}(x,\Gx)]d\Gx}{\cos\pi\Gx-\cos\pi x}, \quad t\in L, \quad 0<x<1.
\label{2.31}
\eeq
Thus, the boundary conditions on the contour $L$ in (\ref{2.9}) are fulfilled, and since $v(e^{i\pi x})=f(x)$, $u(e^{i\pi x})=\Gf(x)$ and
$\R\Gvf_2(t)=u(t)$, $t\in L$, we deduce the following formula for the solution to the characteristic integral equation (\ref{2.1}):
\beq
\Gf(x)=\fr{\sin\pi x}{2}\int_0^1\fr{[a^{2\Gr_1-1}(x,\Gx)+a^{-2\Gr_1+1}(x,\Gx)]f(\Gx)d\Gx}{\cos\pi\Gx-\cos\pi x},
 \quad 0<x<1,
\label{2.34}
\eeq
with the function $a(x,\Gx) $ given by (\ref{2.30}). We wish now to transform the solvability condition (\ref{2.27}). On making the substitution $t=e^{i\pi x}$
and employing the formulas
$$
\fr{\Go_j(t)}{\Go_j(0)}=-ie^{\pi i\Gr_j}\tan^{1-2\Gr_j}\fr{\pi x}{2},
$$
\beq
\cos\fr{\pi\Gd}{2}=\sqrt{\fr{1+|\Gb|}{2}}, \quad 
\sin\fr{\pi\Gd}{2}=\sqrt{\fr{1-|\Gb|}{2}}\sgn\Gb,
\label{2.34.1}
\eeq
we deduce
\beq
\int_0^1\left(\tan^{2\Gr_1-1}\fr{\pi x}{2}+\cot^{2\Gr_1-1}\fr{\pi x}{2}\right)f(x)dx=0.
\label{2.33}
\eeq
We have thus shown that, in the case $-1<\Gb<1$, the integral equation (\ref{2.1}) is solvable 
in the class of bounded at the ends functions if and only if the function $f(x)$
meets the condition (\ref{2.33}). If this condition is satisfied, then 
the solution is unique and given by (\ref{2.34}).
Note that the solution and the solvability condition  derived admit an alternative representation.
On making the substitutions $\Gn=\cos\pi\Gx$, $\Gz=\cos\pi x$, we recast the solution and the condition 
(\ref{2.33}) as
$$
\Gf\left(\fr{\cos^{-1}\Gz}{\pi}\right)=\fr{1}{2\pi}\int_{-1}^1
\left[\left(\fr{1-\Gn}{1+\Gn}\fr{1+\Gz}{1-\Gz}\right)^{\Gr_1-\fr12}+
\left(\fr{1+\Gn}{1-\Gn}\fr{1-\Gz}{1+\Gz}\right)^{\Gr_1-\fr12}\right]\sqrt{\fr{1-\Gz^2}{1-\Gn^2}}
\fr{f\left(\pi^{-1}{\cos^{-1}\Gn}\right)d\Gn}{\Gn-\Gz},
$$
\beq
\int_{-1}^1\left[
\left(\fr{1+\Gz}{1-\Gz}\right)^{\Gr_1-\fr12}+
\left(\fr{1-\Gz}{1+\Gz}\right)^{\Gr_1-\fr12}\right]f\left(\fr{\cos^{-1}\Gz}{\pi}\right)\fr{d\Gz}{\sqrt{1-\Gz^2}}=0,
\label{2.35'}
\eeq
respectively.

 Having defined the exact formula for the solution  let us now show that the function $\Gf(x)$ is not just
bounded as $x\to 0$ and $x\to 1$, but vanishes at these points.
On utilizing
the following formulas  for  the Cauchy integral  when
 $h(\Gn)\in H[-1,1]$ and $0<\R \Ga<1$ {\bf(\ref{mus})}, p.73:
$$
\int_{-1}^1\fr{h(\Gn)d\Gn}{(\Gn+1)^\Ga(\Gn-\Gz)}\sim\pi h(-1)\cot\pi \Ga (\Gz+1)^{-\Ga} +A_1(\Gz), \quad \Gz\to -1^+,
$$
\beq
\int_{-1}^1\fr{h(\Gn)d\Gn}{(1-\Gn)^\Ga(\Gn-\Gz)}\sim -\pi h(1)\cot\pi \Ga (1-\Gz)^{-\Ga} +A_2(\Gz), \quad \Gz\to 1^-,
\label{2.36}
\eeq
we notice that in the representation (\ref{2.35'}) of the function $\Gf(x)$
the bounded terms are canceled, and in the vicinity of the point $\Gz=-1$, the function $\Gf(x)$ behaves as
\beq
\Gf(x)=C_1^-(\Gz+1)^{\Gr_1}+C_2^-(\Gz+1)^{1-\Gr_1}+o((\Gz+1)^{\Gr_0}), 
\quad \Gz\to -1^+,
\label{2.36'}
\eeq
Here, $A_1(\Gz)$ and $A_2(\Gz)$ are functions analytic in a neighborhood
of the points $\Gz=-1$ and $\Gz=1$, respectively,  $C_1^-$, and $C_2^-$ are nonzero constants, and 
\beq
\Gr_0=\min\{\Gr_1,1-\Gr_1\}=1-\Gr_1=\left\{\begin{array}{cc}
\fr12-\fr{\Gd}{2}, &\Gb>0,\\
-\fr{\Gd}{2}, &\Gb<0.\\
\end{array}
\right.
\label{2.37.1}
\eeq
A similar argument is applied to the case when $\Gz\to 1-0$, and we have 
\beq
\Gf(x)\sim C^-(\Gz+1)^{\Gr_0}, \quad \Gz\to -1^+,
\qquad
\Gf(x)\sim  C^+(-\Gz+1)^{\Gr_0}, \quad \Gz\to 1^-.
\label{2.37}
\eeq
with  $C^-$ and $C^+$ being nonzero constants.
Now, since
$$
1+\Gz\sim\fr{\pi^2}{2}(x-1)^2, \quad \Gz\to -1^+, \quad x\to 1^-,
$$
\beq
1-\Gz\sim\fr{\pi^2}{2}x^2, \quad \Gz\to 1^-, \quad x\to 0^+,
\label{2.37.2}
\eeq
we deduce
that the function $\Gf(x)$ vanishes at the points $x=0$ and $x=1$ and
\beq
\Gf(x)\sim C_0^-(1-x)^{2-2\Gr_1}, \quad x\to 1^-,
\qquad
\Gf(x)\sim C_0^+x^{2-2\Gr_1}, \quad x\to 0^+,
\label{2.37.3}
\eeq
where
\beq
2-2\Gr_1=\left\{\begin{array}{cc}
1-\Gd\in\left(\fr12,1\right), &\Gb>0,\\
-\Gd\in\left(0,\fr12\right), &\Gb<0,\\
\end{array}
\right.\quad \Gd=\fr{1}{\pi}\tan^{-1}\fr{\sqrt{1-\Gb^2}}{\Gb}.
\label{2.37.4}
\eeq

\subsection{Case $\Gb=0$}\label{s2.3}

In this particular case the term with  two fixed singularities in the kernel of  equation (\ref{2.1}) vanishes, and the equation becomes
\beq
\fr12\int_0^1
\cot\fr{\pi(\Gx-x)}{2}
\Gf(\Gx)d\Gx=f(x), \quad 0<x<1,
\label{2.38}
\eeq
For the parameters introduced we have the following values:
\beq
\Gb_0=1, \quad \Gl_1=i, \quad \Gl_2=-i, \quad \Gr_1=\fr34, \quad \Gr_2=\fr14, 
\label{2.39}
\eeq
and $\Gd=\pm\fr12$ if $\Gb=0^\pm$.
The solvability condition (\ref{2.33}) 
reduces to the form
\beq
\int_0^1\fr{f(x)}{\sqrt{\sin\pi x}}\cos\left(\fr{\pi x}{2}-\fr{\pi}{4}\right)dx=0. 
\label{2.41}
\eeq
If this condition is satisfied, then according to section \ref{2.2} the solution to equation (\ref{2.38})
reads
\beq
\Gf(x)=\fr12\int_0^1\sqrt{\fr{\sin\pi x}{\sin\pi\Gx}}\fr{f(\Gx)d\Gx}{\sin\fr{\pi}2(x-\Gx)}.
\label{2.42}
\eeq
From here we immediately deduce that the solution vanishes at both ending points,
\beq
\Gf(x)\sim D_0 x^{1/2}, \quad x\to 0^+, \qquad \Gf(x)\sim D_1 (1-x)^{1/2}, \quad x\to 1^-,
\label{2.42'}
\eeq
where $D_0$ and $D_1$ are nonzero constants.

We show now that the results found are  consistent with the ones  recovered from the classical theory  {\bf(\ref{gak})}.
First, by making the substitutions $t=e^{i\pi x}$, $\tau=e^{i\pi\Gx}$ and denoting $\tilde\Gf(t)=\Gf(x)$, $\tilde f(t)=f(x)$, we
rewrite equation (\ref{2.38}) as
\beq
\fr{1}{\pi}\int_L\fr{\tilde\Gf(\tau)d\tau}{\tau-t}=\tilde f(t)+\tilde\Gf_0,\quad t\in L,
\label{2.43}
\eeq
where
\beq
\tilde\Gf_0=\fr{1}{2\pi}\int_L\fr{\tilde\Gf(t)dt}{t},
\label{2.44}
\eeq
and, as before, $L=\{z\in {\Bbb C}: |z|=1, \I z>0\}$ with the starting point $z=1$. In the class of functions bounded at the endpoints of the contour $L$ the solution to equation (\ref{2.43}) does not exist unless
the condition 
\beq
\int_L\fr{\tilde f(t)+\tilde\Gf_0}
{\sqrt{t^2-1}}dt=0
\label{2.45}
\eeq
is fulfilled. Then the solution is unique and 
has the form {\bf(\ref{gak})}
\beq
\tilde\Gf(t)=-\fr{1}{\pi}\sqrt{t^2-1}\int_L\fr{\tilde f(\tau)+\tilde\Gf_0}{\sqrt{\tau^2-1}}\fr{d\tau}{\tau-t}.
\label{2.46}
\eeq
Here, the branch of the square root in the plane cut along the contour $L$ 
is fixed by the condition $\sqrt{z^2-1}\sim z$, $z\to\infty$. Note that the lower 
side of the cut $|z|=1-0$, $\arg z\in[0,\pi]$, is identified as the contour $L$ itself.
Employing the relation
\beq
\int_L\fr{d\tau}
{\sqrt{\tau^2-1}(\tau-t)}=0, \quad t\in L,
\label{2.47}
\eeq
and coming back to the original variables and functions we transform formula (\ref{2.46}) into the form (\ref{2.42}).
Analyze now the solvability condition (\ref{2.45}). Owing to (\ref{2.44}) and (\ref{2.46})
we thereby represent the constant $\tilde\Gf_0$ as the double integral
\beq
\tilde\Gf_0=-\fr{1}{2\pi^2}\int_L\fr{\sqrt{t^2-1}dt}{t}\int_L\fr{\tilde f(\tau)d\tau}{\sqrt{\tau^2-1}(\tau-t)}.
\label{2.48}
\eeq
To convert the double integral into a single one, we introduce the complex potential
\beq
\GO(z)=\fr{1}{2\pi i}\int_L\fr{\sqrt{\tau^2-1}d\tau}{\tau(\tau-z)}. \quad z\ne 0.
\label{2.49}
\eeq
Applying the theory of residues and taking into account that for the branch fixed,
$\sqrt{z^2-1}|_{z=0}=-i$, we discover
\beq
\GO(z)=-\fr12+\fr{i}{2z}+\fr{\sqrt{z^2-1}}{2z}.
\label{2.50}
\eeq
To determine the principal value $\GO(t)$ of the integral (\ref{2.49}), we apply the Sokhotski-Plemelj formulas and obtain
\beq
\fr{1}{2\pi i}\int_L\fr{\sqrt{\tau^2-1}d\tau}{\tau(\tau-t)}=-\fr12+\fr{i}{2t}.
\label{2.51}
\eeq
Therefore
\beq
\tilde\Gf_0=-\fr{1}{2\pi}\int_L\fr{\tilde f(\tau)}{\sqrt{\tau^2-1}}\left(i+\fr{1}{\tau}\right)d\tau.
\label{2.52}
\eeq
Next, by substituting this expression in (\ref{2.45}) and making use of the integral 
\beq
\int_L\fr{d\tau}{\sqrt{\tau^2-1}}=-\pi i
\label{2.53}
\eeq
we discover that $\tilde\Gf_0=0$ that is the solvability condition becomes
\beq
\int_L\fr{\tilde f(\tau)}{\sqrt{\tau^2-1}}\left(i+\fr{1}{\tau}\right)d\tau=0.
\label{2.53'}
\eeq
We finally come back to the variable $x$ and the function $f(x)$ and  
reduce the condition (\ref{2.53'}) to the desired form (\ref{2.41}).

\subsection{Case $\Gb>1$}\label{s2.4}

Now we assume that $\Gb>1$. In this case  the eigenvalues of the matrix $G_0$,
$\Gl_j=-\Gb-(-1)^j\sqrt{\Gb^2-1}$, $j=1,2$, are real and negative. 
The entries of the matrix of transformation $T$ are also real,
\beq
T=\left(\begin{array}{cc}
1 & 1\\
\sqrt{-\Gb_0} & -\sqrt{-\Gb_0} \\
\end{array}\right), 
\quad \Gb_0=\fr{1+\Gb}{1-\Gb}<0.
\label{2.54}
\eeq
The matrix of factorization $X(s)$ has the same form as (\ref{2.22.3}).
However the parameters $\Gr_1$ and $\Gr_2$ are not real anymore,
\beq
\Gr_1=\fr12-i\Gve, \quad \Gr_2=\fr12+i\Gve,
\label{2.55}
\eeq
where
\beq
\Gve=\fr{1}{2\pi}\log(\Gb+\sqrt{\Gb^2-1})>0.
\label{2.56}
\eeq
Following the scheme of section \ref{s2.2} we derive the solution of the Hilbert problem (\ref{2.9}) in the form
$$
\Gvf_1(z)=
\fr{z^2-1}{2\pi \sqrt{-\Gb_0}}\int_L
\left[\fr{\Go_1(t)}{\Go_1(z)}-\fr{\Go_2(t)}{\Go_2(z)}\right]\fr{v(t)dt}{(t-z)(1-tz)},
$$
\beq
\Gvf_2(z)=\fr{z^2-1}{2\pi}\int_L\left[\fr{\Go_1(t)}{\Go_1(z)}+\fr{\Go_2(t)}{\Go_2(z)}\right]\fr{v(t)dt}{(t-z)(1-tz)},
\label{2.57}
\eeq
where 
\beq
\Go_j(t)=ie^{-\pi i\Gr_j}\tan^{1-2\Gr_j}\fr12\pi x,\quad  j=1,2, \quad
1-2\Gr_1=2i\Gve, \quad  1-2\Gr_2=-2i\Gve.
\label{2.57'}
\eeq

The condition (\ref{2.10}) that guarantees that the Hilbert problem (\ref{2.9}) 
is equivalent to the integral equation (\ref{2.4}) reads
\beq
\int_0^1\left(\tan^{2i\Gve}\fr{\pi x}{2}+\cot^{2i\Gve}\fr{\pi x}{2}\right) v(e^{i\pi x})dx=0.
\label{2.58}
\eeq
Here we used the relations
\beq
e^{\pi\Gve}=\fr{1}{\sqrt{\Gb-\sqrt{\Gb^2-1}}}, \quad
\left(\fr{1}{\sqrt{-\Gb_0}}\mp 1\right)e^{\pm\pi\Gve}=\mp\sqrt{\fr{2}{\Gb+1}}.
\label{2.59}
\eeq
We assume further that the condition (\ref{2.58}) is satisfied. Our task now is to
find the solution to the integral equation (\ref{2.1}). The Sokhotski-Plemelj formulas 
applied to (\ref{2.57}) yield 
$$
\Gvf_1(e^{i\pi x})=
\fr{i\sin\pi x}{\sqrt{-\Gb_0}}\int_0^1\fr{v(e^{i\pi \Gx})}{\cos\pi\Gx-\cos\pi x}
\sin\left(2\Gve\log a(x,\Gx)\right)d\Gx,
$$
\beq
\Gvf_2(e^{i\pi x})=-iv(e^{i\pi x})+\sin\pi x\int_0^1\fr{ v(e^{i\pi \Gx})}{\cos\pi\Gx-\cos\pi x}
\cos\left(2\Gve\log a(x,\Gx)\right)d\Gx,
\label{2.60}
\eeq
where $a(x,\Gx)$ is given by (\ref{2.30}).
These formulas enable us to verify the boundary conditions 
of the Hilbert problem ({\ref{2.9}) on the contour $L$
\beq
\R\Gvf_1(t)=0,\quad \I \Gvf_2(t)=-v(t), \quad t\in L,
\label{2.61}
\eeq
and also to derive the solution to the integral equation (\ref{2.4}), $u(t)=\R\Gvf_2(t).$
On putting $f(\Gx)=v(e^{i\pi\Gx})$ and  $\Gf(x)=u(e^{i\pi x})$ we have from (\ref{2.60})
\beq
\Gf(x)=\int_0^1\fr{\sin\pi x f(\Gx)}{\cos\pi\Gx-\cos\pi x}
\cos\left(2\Gve\log a(x,\Gx)\right)d\Gx, \quad 0<x<1.
\label{2.62}
\eeq
We assert that the integral equation (\ref{2.1}) has at most one solution  in the class of functions 
bounded at the ends. If the condition 
\beq
\int_0^1\cos\left(2\Gve\log\tan\fr{\pi x}{2}\right)f(x)dx=0
\label{2.63}
\eeq
is satisfied, then the solution exists and is given by (\ref{2.62}).
Alternatively,  applying the map $\Gz=\cos\pi x$, $\Gn=\cos\pi\Gx$ we can write this condition and the solution as
$$
\int_{-1}^1\left[
\left(\fr{1+\Gz}{1-\Gz}\right)^{i\Gve}+
\left(\fr{1+\Gz}{1-\Gz}\right)^{-i\Gve}\right]f\left(\fr{\cos^{-1}\Gz}{\pi}\right)\fr{d\Gz}{\sqrt{1-\Gz^2}}=0,
$$
\beq
\Gf\left(\fr{\cos^{-1}\Gz}{\pi}\right)=\fr{1}{2\pi}\int_{-1}^1
\left[\left(\fr{1-\Gn}{1+\Gn}\fr{1+\Gz}{1-\Gz}\right)^{i\Gve}+
\left(\fr{1-\Gn}{1+\Gn}\fr{1+\Gz}{1-\Gz}\right)^{-i\Gve}\right]\sqrt{\fr{1-\Gz^2}{1-\Gn^2}}
\fr{f\left(\pi^{-1}{\cos^{-1}\Gn}\right)d\Gn}{\Gn-\Gz},
\label{2.63'}
\eeq
respectively. 
To discover the asymptotics of the solution at the endpoints, we employ formulas (\ref{2.36})
as we did in the case $-1<\Gb<1$. At the point $\Gz=-1$, we have
\beq
\Gf(x)\sim C_0(1+\Gz)^{1/2+i\Gve}+\ov{C_0}(1+\Gz)^{1/2-i\Gve}, \quad \Gz\to-1^+,
\label{2.63.0}
\eeq
and then, due to (\ref{2.37.2}), derive
\beq
\Gf(x)\sim  (1-x)[C_{10}\cos (2\Gve\log (1-x))+C_{11}\sin (2\Gve\log (1-x))], \quad x\to 1^-,
\label{2.63''}
\eeq
Similarly, at the second endpoint $x=0$,
\beq
\Gf(x)\sim  x[C_{00}\cos (2\Gve\log x)+C_{01}\sin (2\Gve\log x)], \quad x\to 0^+.
\label{2.63'''}
\eeq
Here, $C_{jm}$ ($j,m=0,1$) are real nonzero constants.

\subsection{Case $\Gb<-1$}\label{s2.5}

In this case, the eigenvalues $\Gl_1$ and $\Gl_2$ of the matrix $G_0$ are positive,
$\Gl_j=-\Gb-(-1)^j\sqrt{\Gb^2-1}$, $j=1,2$,    and 
the transformation matrix $T$ 
has the form
\beq
T=\left(\begin{array}{cc}
1 & 1\\
-\sqrt{-\Gb_0} & \sqrt{-\Gb_0} \\
\end{array}\right), 
\quad \Gb_0=\fr{1+\Gb}{1-\Gb}<0.
\label{2.65}
\eeq
There are two possibilities for the choice of the parameters $\Gr_1$ and $\Gr_2$, (i) $ \Gr_1=i\Gve,  \Gr_2=-i\Gve,$
and (ii) $\Gr_1=1+i\Gve,  \Gr_2=1-i\Gve.$ Here,
\beq
 \Gve=\fr{1}{2\pi}\log(-\Gb+\sqrt{\Gb^2-1})>0.
\label{2.64}
\eeq
Note that the pairs $ \Gr_1=i\Gve$,  $\Gr_2=1-i\Gve$ and $ \Gr_1=1+i\Gve$,  $\Gr_2=-i\Gve$
lead to the functions $\Gvf_1(z)$ and $\Gvf_2(z)$ which do not satisfy the boundary 
conditions (\ref{2.61}).

For the  pair  $ \Gr_1=i\Gve$,  $\Gr_2=-i\Gve$, 
similarly to sections \ref{2.2} and \ref{2.4}, we derive the solution to the Hilbert problem
(\ref{2.9}) in the form
$$
\Gvf_1(z)=
-\fr{z^2-1}{2\pi \sqrt{-\Gb_0}}\int_L
\left[\fr{\Go_1(t)}{\Go_1(z)}-\fr{\Go_2(t)}{\Go_2(z)}\right]\fr{v(t)dt}{(t-z)(1-tz)},
$$
\beq
\Gvf_2(z)=\fr{z^2-1}{2\pi}\int_L\left[\fr{\Go_1(t)}{\Go_1(z)}+\fr{\Go_2(t)}{\Go_2(z)}\right]\fr{v(t)dt}{(t-z)(1-tz)},
\label{2.66}
\eeq
where 
\beq
\Go_1(t)=ie^{\pi\Gve}\tan^{1-2i\Gve}\fr{\pi x}{2}, \quad
\Go_2(t)=ie^{-\pi\Gve}\tan^{1+2i\Gve}\fr{\pi x}{2}.
\label{2.67}
\eeq
As before, by applying the Sokhotski-Plemelj formulas  to the singular integrals (\ref{2.66}) we determine
$$
\Gvf_1(e^{i\pi x})=
\fr{i\sin\pi x}{\sqrt{-\Gb_0}}\int_0^1\fr{v(e^{i\pi \Gx})a(x,\Gx)}{\cos\pi\Gx-\cos\pi x}
\sin\left(2\Gve\log a(x,\Gx)\right)d\Gx,
$$
\beq
\Gvf_2(e^{i\pi x})=-iv(e^{i\pi x})+\sin\pi x\int_0^1\fr{ v(e^{i\pi \Gx})a(x,\Gx)}{\cos\pi\Gx-\cos\pi x}
\cos\left(2\Gve\log a(x,\Gx)\right)d\Gx,
\label{2.68}
\eeq
with $a(x,\Gx)$ being the function given by (\ref{2.30}).
The functions (\ref{2.68}) satisfy the boundary conditions on $L$ in 
(\ref{2.9}), $\R\Gvf_1(t)=0$, $\I \Gvf_2(t)=-v(t)$,
and yield the solution to the integral equation (\ref{2.1}), $\Gf(x)=\R\Gvf_2(e^{i\pi x})$ with
$v(e^{i\pi \Gx})=f(\Gx)$,
\beq
\Gf(x)=\int_0^1\fr{\sin\pi x a(x,\Gx)
\cos\left(2\Gve\log a(x,\Gx)\right)}{\cos\pi\Gx-\cos\pi x}f(\Gx)d\Gx, \quad 0<x<1.
\label{2.69}
\eeq
For this function to be the solution to 
the integral equation (\ref{2.1}) in the class of functions chosen, it is necessary and sufficient that
the function $f(x)$ meets the condition (\ref{2.10}) or, equivalently,
\beq
\int_0^1\tan\fr{\pi x}{2}\cos\left(2\Gve\log\tan\fr{\pi x}{2}\right)f(x)dx=0.
\label{2.70}
\eeq
This condition can be also represented in the form
\beq
\int_{-1}^1\left[
\left(\fr{1+\Gz}{1-\Gz}\right)^{i\Gve}+
\left(\fr{1+\Gz}{1-\Gz}\right)^{-i\Gve}\right]f\left(\fr{\cos^{-1}\Gz}{\pi}\right)\fr{d\Gz}{1+\Gz}=0,
\label{2.70'}
\eeq
and the corresponding solution to the integral equation becomes 
\beq
\Gf\left(\fr{\cos^{-1}\Gz}{\pi}\right)=\fr{1}{2\pi}\int_{-1}^1
\left[\left(\fr{1-\Gn}{1+\Gn}\fr{1+\Gz}{1-\Gz}\right)^{i\Gve}+
\left(\fr{1-\Gn}{1+\Gn}\fr{1+\Gz}{1-\Gz}\right)^{-i\Gve}\right]
\fr{f\left(\pi^{-1}{\cos^{-1}\Gn}\right)}{\Gn-\Gz} \fr{(1+\Gz)d\Gn}{1+\Gn}.
\label{2.70''}
\eeq
Its analysis shows that the solution vanishes at the point $\Gz=-1$, while at the second end, $\Gz=1$,
it is bounded and oscillates. 
Indeed, in view of the relations (\ref{2.36}) applied to (\ref{2.70''}) we have
$$
\Gf\left(\fr{\cos^{-1}\Gz}{\pi}\right)\sim D_1(1+\Gz)^{1+i\Gve}+ \ov{D_1}(1+\Gz)^{1-i\Gve},
\quad \Gz\to -1^+,
$$
\beq\Gf\left(\fr{\cos^{-1}\Gz}{\pi}\right)\sim D_0(1-\Gz)^{i\Gve}+ \ov{D_0}(1-\Gz)^{-i\Gve},
\quad \Gz\to 1^-.
\label{2.70'''}
\eeq
for some nonzero complex constants $D_0$ and $D_1$.
On returning to the original variables we 
deduce 
$$
\Gf(x)\sim  D_{00}\cos (2\Gve\log x)+D_{01}\sin (2\Gve\log x), \quad x\to 0^+,
$$
\beq
\Gf(x)\sim  (1-x)^2[D_{10}\cos (2\Gve\log (1-x))+D_{11}\sin (2\Gve\log (1-x))], \quad x\to 1^-,
\label{2.70.0}
\eeq
where $D_{jm}$ ($j,m=0,1$) are real nonzero constants.

For the second choice of the parameters $\Gr_1$ and $\Gr_2$, $\Gr_1=1+i\Gve$ and $\Gr_2=1-i\Gve$, in the case 
$\Gb<-1$, we again derive the solution to the Hilbert problem (\ref{2.9}), (\ref{2.10}), take 
the real part of the second potential $\Gvf_2(z)$, put $v(e^{i\pi x})=f(x)$ and obtain the solvability
condition of the equation (\ref{2.1})
\beq
\int_0^1\cot\fr{\pi x}{2}\cos\left(2\Gve\log\tan\fr{\pi x}{2}\right)f(x)dx=0
\label{2.70.1}
\eeq  
and the solution to the integral equation
\beq
\Gf(x)=\int_0^1\fr{\sin\pi x 
\cos\left(2\Gve\log a(x,\Gx)\right)}{(\cos\pi\Gx-\cos\pi x)a(x,\Gx)}f(\Gx)d\Gx, \quad 0<x<1.
\label{2.70.2}
\eeq
On making the substitutions $\Gn=\cos\pi\Gx$, $\Gz=\cos\pi x$, we have another form of the solvability 
condition and the solution to the integral equation. It coincides with formulas (\ref{2.70'}) and (\ref{2.70''})
if we replace there $d\Gz/(1+\Gz)$ and $(1+\Gz)d\Gn/(1+\Gn)$ by  $d\Gz/(1-\Gz)$ and 
$(1-\Gz)d\Gn/(1-\Gn)$, respectively.
The solution (\ref{2.70.2}) associated with the parameters $\Gr_1=1+i\Gve$ and $\Gr_2=1-i\Gve$
and the condition (\ref{2.70.1}) vanishes at the point $x=0$ and is bounded and oscillates at the point $x=1$,
$$
\Gf(x)\sim x^2 [E_{00}\cos (2\Gve\log x)+E_{01}\sin (2\Gve\log x)], \quad x\to 0^+,
$$
\beq
\Gf(x)\sim  E_{10}\cos (2\Gve\log (1-x))+E_{11}\sin (2\Gve\log (1-x)), \quad x\to 1^-,
\label{2.70.3}
\eeq
where $E_{jm}$ ($j,m=0,1$) are real nonzero constants.

\subsection{Solution to the characteristic integral equation}\label{s2.6}

We now summarize the results.
\vspace{.1in}

Theorem 2.2. Let $f(x)\in H[0,1]$ and $\Gb$ be a real parameter. Denote
$$
\Gve=\left\{
\begin{array}{cc}
\Gd, & \Gb\in(0,1),\\
1+\Gd, & \Gb\in(-1,0),\\
(2\pi)^{-1}\log(|\Gb|+\sqrt{\Gb^2-1}), & |\Gb|>1,
\end{array}
\right.\quad \Gd=\fr{1}{\pi}\tan^{-1}\fr{\sqrt{1-\Gb^2}}{\Gb},
$$
\beq
a(x,\Gx)=\tan\fr{\pi \Gx}{2}\cot\fr{\pi x}{2}.
\label{2.71}
\eeq
In the class of functions H\"older-continuous
everywhere in the interval $(0,1)$ and bounded at the endpoints,  generally, the singular integral equation
\beq
\CS[\Gf](x)\equiv \int_0^1S(x,\Gx)\Gf(\Gx)d\Gx=f(x), \quad 0<x<1,
\label{2.71'}
\eeq
with the kernel $S(x,\Gx)=\fr12\cot\fr{\pi}2(\Gx-x)+\fr{\Gb}2\cot\fr{\pi}2(\Gx+x)$
does not have a
solution. If $|\Gb|<1$ or $\Gb>1$, the integral equation becomes solvable if and only if the function $f(x)$ meets the 
condition
\beq
\int_0^1 V(x)f(x)dx=0,
\label{2.76'}
\eeq
where
\beq
V(x)=\left\{
\begin{array}{cc}
(\sin\pi x)^{-1/2}\cos\left(\fr{\pi x}{2}-\fr{\pi}{4}\right), & \Gb=0,\\
\tan^{\Gve}\fr{\pi x}{2}+\cot^{\Gve}\fr{\pi x}{2}, & |\Gb|<1, \;\Gb\ne 0,\\
\cos\left(2\Gve\log\tan\fr{\pi x}{2}\right), & \Gb>1.\\
\end{array}
\right.
\label{2.76''}
\eeq
If the condition of solvability is fulfilled, then the inverse operator $S^{-1}$ exists, and  the solution is unique
and given by $\Gf(x)=\CS^{-1}[f](x)$, where
$$
\CS^{-1}[f](x)=\fr12\int_0^1\sqrt{\fr{\sin\pi x}{\sin\pi\Gx}}\fr{f(\Gx)d\Gx}{\sin\fr12\pi(x-\Gx)}, \quad \Gb=0,
$$
$$
\CS^{-1}[f](x)=\fr{1}{2}\int_0^1\left[a^{\Gve}(x,\Gx)+a^{-\Gve}(x,\Gx)\right]
\fr{\sin\pi x f(\Gx)d\Gx}{\cos\pi\Gx-\cos\pi x},\quad |\Gb|<1, \quad \Gb\ne 0,
$$
\beq
\CS^{-1}[f](x)=\int_0^1
\cos\left(2\Gve\log a(x,\Gx)
\right)\fr{\sin\pi x f(\Gx)d\Gx}{\cos\pi\Gx-\cos\pi x},
\quad \Gb>1.
\label{2.77}
\eeq
The solution vanishes at the points $x=0$ and $x=1$, and its asymptotics is described by (\ref{2.42'}) if
$\Gb=0$, (\ref{2.37.3}) in the case
$|\Gb|<1$ and (\ref{2.63''}), (\ref{2.63'''}) if
$\Gb>1$.

In the case $\Gb<-1$,  in the class of functions chosen, the integral equation is solvable if and only if the function $f(x)$ satisfies one of the following two conditions:
\beq 
\int_0^1\left(\tan\fr{\pi x}{2}\right)^{\pm 1}\cos\left(2\Gve\log\tan\fr{\pi x}{2}\right)f(x)dx=0.
\label{2.77.0}
\eeq 
Then the solution  is unique and given by $\Gf(x)=\CS^{-1}[f](x)$, where
\beq
\CS^{-1}[f](x)=\int_0^1\fr{\sin\pi x [a(x,\Gx)]^{\pm 1}
\cos\left(2\Gve\log a(x,\Gx)\right)}{\cos\pi\Gx-\cos\pi x}f(\Gx)d\Gx, \quad 0<x<1.
\label{2.77.1}
\eeq
It is bounded and oscillates at the point $p^\pm$ and vanishes 
at the point $p^\mp$, where $p^+=0$ and $p^-=1$. The asymptotics of the solution is described
by (\ref{2.70.0}) and (\ref{2.70.3}).

\vspace{.1in}

The two cases left, $\Gb=1$ and $\Gb=-1$, can be easily treated by reducing the integral equation
to an equation solvable by the Hilbert inversion formula (see appendix B).

\setcounter{equation}{0}

\section{Complete singular integral equation}\label{s3}

In this section, we will develop an algorithm for the solution of the complete
singular integral equation
\beq
\int_0^1[S(x,\Gx)+K(x,\Gx)]\Gf(\Gx)d\Gx=-F(x)+C, \quad 0<x<1,
\label{3.1}
\eeq
where  $S(x,\Gx)$ is the singular kernel introduced in (\ref{2.2}), $K(x,\Gx)$ is a regular kernel,
$F(x)\in H[0,1]$, and $C$ is an unknown constant. We seek the solution to this equation, $\Gf(x)$,
 in the class of H\"older functions bounded at the endpoints. 

\subsection{Regularization of the complete equation}

First we regularize equation (\ref{3.1}) and reduce it to a Fredholm integral equation. 
Since the inverse operator   $\CS^{-1}$ has been constructed, it is reasonable to apply the 
Carleman-Vekua regularization procedure  {\bf(\ref{gak})}. 
This brings us to the equation
\beq
\Gf(x)+\CS^{-1}\CK[\Gf](x)=\CS^{-1}[-F+C](x), \quad 0<x<1,
\label{3.1.1}
\eeq
provided the constant $C$ is selected to be
\beq
C=\left(\int_0^1 V(x)dx\right)^{-1}\int_0^1 V(x)\left[F(x)+\CK[\Gf](x)\right]dx.
\label{3.1.2}
\eeq
Here, $\CK$ is the Fredholm operator with the kernel $K(x,\Gx)$,   and $V(x)$ is the function defined by (\ref{2.76''}).
 Motivated by applications
 to fracture mechanics  
 we restrict ourselves to considering
 the case $|\Gb|<1$ and exclude the trivial case $\Gb=0$. Note that a similar procedure can be worked out when $|\Gb|>1$.
 So, here and further we suppose $0<|\Gb|<1$.  It turns out (see Appendix C) that 
 the first integral in (\ref{3.1.2}) can be evaluated explicitly and the solvability condition
 becomes
 \beq
 C=\fr{\sin\pi \Gr_1}{2}\int_0^1\left(\tan^{2\Gr_1-1}\fr{\pi x}{2}+\cot^{2\Gr_1-1}\fr{\pi x}{2}\right)\left[F(x)+\CK[\Gf](x)\right]dx.
 \label{3.1.3}
 \eeq
 If this condition holds, then the integral equation (\ref{3.1}) is equivalent to the following Fredholm
 integral equation:
 \beq
 \Gf(x)+\int_0^1 \tilde K(x,\Gx)\Gf(\Gx)d\Gx=-\tilde F(x), \quad 0<x<1,
 \label{3.1.4}
 \eeq
 where $\tilde K(x,\Gx)$ and $\tilde F(x)$ are bounded in the sets $0\le x,\Gx\le 1$ and $0\le x\le 1$, respectively, and have the form
 $$
 \tilde K(x,\Gx)=\fr{\sin\pi x}{2}\int_0^1\fr{[a^{2\Gr_1-1}(x,\tau)+a^{-2\Gr_1+1}(x,\tau)]K(\tau,\Gx)d\tau}
 {\cos\pi\tau-\cos\pi x},
$$
\beq
\tilde F(x)= \fr{\sin\pi x}{2}\int_0^1\fr{[a^{2\Gr_1-1}(x,\tau)+a^{-2\Gr_1+1}(x,\tau)]F(\tau)d\tau}
 {\cos\pi\tau-\cos\pi x}.
 \label{3.1.5}
 \eeq
 Here, we used the relation $\CS^{-1}[1]=0$  to be proved in Section \ref{ss3.1}.
As in the case of the complete singular integral equation with the Cauchy kernel, for numerical purposes, 
it is preferable to deal with the singular equation (\ref{3.1}) directly and bypass its regularization.
In the next sections we generalize  the method of orthogonal  
polynomials efficient for equations with the Cauchy kernel to the case of equation
(\ref{3.1}).

\subsection{Spectral relation for the operator $\CS$}\label{ss3.1} 

The heart of the numerical method to be proposed is the derivation of the solutions $\Gf_j(x)$,   
$j=0,1,\ldots$,
of the characteristic equation $\CS[\Gf_j](x)=f_j(x)$ with the right-hand side chosen to be $f_j(x)=C_j-\cos[(j+1)\pi x]$,
\beq
\fr12\int_0^1
\left[\cot\fr{\pi(\Gx-x)}{2}+\Gb\cot\fr{\pi(\Gx+x)}{2}
\right]\Gf_j(\Gx)d\Gx=C_j-\cos[(j+1)\pi x], \quad 0<x<1, \quad j=0,1,\ldots,
 \label{3.4}
 \eeq
where $C_j$ are constants to be fixed.
The solutions are sought in the class of functions vanishing at the endpoints and meeting the condition  (\ref{2.33}) that is 
$\CS: L_w^2(0,1)\cap H^\circ(0,1)\to L^2(0,1)\cap H[0,1]$, 
where $L_w^2(0,1)$
is the weight Hilbert space,  $w(x)=\sin^{2(\Gr_1-1)} \pi x$, and $H^\circ$
is the class of H\" older functions  meeting the condition  (\ref{2.33}).
The functions $\Gf_j(x)$ admit the following explicit expressions in terms of some
trigonometric polynomials:
\beq
\Gf_j(x)=\cos^{2\Gr_1}\fr{\pi x}{2}\sin^{2(1-\Gr_1)}\fr{\pi x}{2} q_j^{(\Gr_1)}(x)+
\cos^{2(1-\Gr_1)}\fr{\pi x}{2}\sin^{2\Gr_1}\fr{\pi x}{2} q_j^{(1-\Gr_1)}(x),
 \label{3.5}
 \eeq
 where
 $$
 q_j^{(\Ga)}(x)=\sum_{\nu=0}^jc_{j\nu}^{(\Ga)}\sin^{2\nu}\fr{\pi x}{2},
 $$
 \beq
 c_{j\nu}^{(\Ga)}=
 \fr{1}{2\sin\pi\Ga}\sum_{m=\nu+1}^{j+1}\fr{(-j-1)_m(j+1)_m(\Ga)_{m-1-\nu}}{(1/2)_m m!(m-1-\nu)!},
 \label{3.6}
 \eeq
and $(\cdot)_m$ is the factorial symbol. On making the substitution $\Gz=\cos\pi x$ it is possible to
rewrite the spectral relation (\ref{3.4}) in the form
$$
\fr1{2\pi}\int_{-1}^1 \left[\fr{\sqrt{(1+\Gn)(1+\Gz)}+\sqrt{(1-\Gn)(1-\Gz)}}
{\sqrt{(1-\Gn)(1+\Gz)}-\sqrt{(1+\Gn)(1-\Gz)}}+\Gb
\fr{\sqrt{(1+\Gn)(1+\Gz)}-\sqrt{(1-\Gn)(1-\Gz)}}
{\sqrt{(1-\Gn)(1+\Gz)}+\sqrt{(1+\Gn)(1-\Gz)}}
\right]\fr{\Gy_j(\Gn)d\Gn}{\sqrt{1-\Gn^2}}
$$
\beq
=C_j-T_{j+1}(\Gz), \quad -1<\Gz<1,
\label{3.7}
\eeq   
where the functions $\Gy_j(\Gz)=\Gf_j(\fr1{\pi}\cos^{-1}\Gz)$ are given by
$$
\Gy_j(\Gz)=\fr12\left[(1+\Gz)^{\Gr_1}(1-\Gz)^{1-\Gr_1}p_j^{(\Gr_1)}(\Gz)+
(1+\Gz)^{1-\Gr_1}(1-\Gz)^{\Gr_1}p_j^{(1-\Gr_1)}(\Gz)\right],
$$
\beq
 p_j^{(\Ga)}(\Gz)=
 \sum_{\nu=0}^jc_{j\nu}
 \left(\fr{1-\Gz}{2}\right)^{\nu}.
\label{3.8}
\eeq
To prove the correctness of these formulas, first we must  satisfy the solvability condition (\ref{2.76'}). It fixes the constants $C_j$
\beq
C_j=\fr{M_{j+1}}{M_0},
 \label{3.9}
 \eeq
 where
 \beq
 M_j=\int_0^1
 \left(\tan^{2\Gr_1-1}\fr{\pi x}{2}+\cot^{2\Gr_1-1}\fr{\pi x}{2}\right)
 T_j(\cos\pi x)dx, \quad j=0,1,\ldots.
 \label{3.10}
 \eeq
 These integrals evaluated explicitly in Appendix C have the form
$$
M_0=2\csc\pi\Gr_1,\quad 
M_{2m-1}=0, 
$$
\beq
 M_{2m}=\fr{2}{\sin\pi\Gr_1}\sum_{j=0}^{2m}\fr{(-2m)_j(2m)_j(\Gr_1)_j}{(1/2)_j(j!)^2},
 \quad m=1,2,\ldots.
 \label{3.15}
\eeq

After we  have computed the constants $C_j$ we wish  to simplify the expression
(\ref{2.34}) ($|\Gb|<1$) for the solution of the integral equation (\ref{3.4}) and therefore verify the representation (\ref{3.5}). On making
the substitutions $\Gz=\cos\pi x$ and $\Gn=\cos\pi\Gx$ and using the first formula in (\ref{2.35'}) we discover
\beq
\Gf_j(x)=\fr12\left[(1+\Gz)^{\Gr_1}(1-\Gz)^{1-\Gr_1}I_j^{(\Gr_1)}(\Gz)+
(1+\Gz)^{1-\Gr_1}(1-\Gz)^{\Gr_1}I_j^{(1-\Gr_1)}(\Gz)\right],
\label{3.16}
\eeq
where
$$
I_j^{(\Ga)}(\Gz)=C_jJ_0^{(\Ga)}(\Gz)-J_{j+1}^{(\Ga)}(\Gz),
$$
\beq
J_{j}^{(\Ga)}(\Gz)=\fr{1}{\pi}\int_{-1}^1(1-\Gn)^{\Ga-1}(1+\Gn)^{-\Ga}\fr{T_j(\Gn)d\Gn}{\Gn-\Gz},
\quad 0<\Ga<1, \quad j=0,1,\ldots.
\label{3.17}
\eeq
The method we apply to compute the latter integral is built upon the convolution
theorem and the theory of residues. Recast the integral $J_{j}^{(\Ga)}(\Gz)$ as
\beq
J_{j}^{(\Ga)}(\Gz)=\int_0^\infty h_1(\tau)h_2\left(\fr{t}{\tau}\right)\fr{d\tau}{\tau},\quad
-1<\Gz<1, \quad 0<t<1,
\label{3.18}
\eeq
where $t=(1-\Gz)/2$,
\beq
h_1(\tau)=\left\{\begin{array}{cc}
\tau^{\Ga-1}(1-\tau)^{-\Ga}T_j(1-2\tau), & 0<\tau<1,\\
0, & \tau>1,\\
\end{array}
\right.\quad
h_2(\tau)=-\fr{1}{2\pi(1-\tau)},
\label{3.19}
\eeq
and apply the convolution theorem. It reads
\beq
J_{j}^{(\Ga)}(\Gz)=\fr{1}{2\pi i}\int_{c-i\infty}^{c+i\infty} H_1(s) H_2(s) t^{-s}ds.
\label{3.20}
\eeq
Here, $H_1$ and $H_2$ are the Mellin transforms of the functions
$h_1$ and $h_2$, respectively, and $c\in(1-\Ga,1)$ ($0<\Ga<1$).
 Formula (\ref{3.13}) enables us to find a series
representation of $H_1(s)$,
\beq
H_1(s)=\fr{\GG(s+\Ga-1)\GG(1-\Ga)}{\GG(s)}{}_3 F_2(-j,j,s+\Ga-1;1/2,s; 1), \quad \R s>1-\Ga,
\label{3.21}
\eeq
and since
\beq
H_2(s)=-\fr12\cot\pi s, \quad 0<\R s<1,
\label{3.22}
\eeq
the integral $J_{j}^{(\Ga)}(\Gz)$ becomes 
\beq
J_{j}^{(\Ga)}(\Gz)=-\fr{\GG(1-\Ga)}{2\pi}\sum_{m=0}^j\fr{(-1)^m(-j)_m (j)_m}{(1/2)_m m!}L_m(t),
\label{3.23}
\eeq
where 
\beq
L_m(t)=\fr{1}{2\pi i} \int_{c-i\infty}^{c+i\infty} \cos\pi s \GG(1-s-m)\GG(s+\Ga-1+m) t^{-s} ds, \quad 0<t<1.
\label{3.24}
\eeq
Next, by employing the theory of residues we compute the integral $L_m(t)$ and deduce the series representation of the function
$J_{j}^{(\Ga)}(\Gz)$
\beq
J_{j}^{(\Ga)}(\Gz)=-\fr{1}{2\sin\pi\Ga}\sum_{m=0}^j\fr{(-j)_m (j)_m t^{m-1}}{(1/2)_m m!}\left[\sum_{k=0}^{m-1}
\fr{(\Ga)_k t^{-k}}{k!}-\cos\pi\Ga\left(\fr{t}{1-t}\right)^\Ga\right].
\label{3.24'}
\eeq
Note that the sum $\sum_{k=0}^{m-1}
(\Ga)_k t^{-k}(k!)^{-1}$ is equal to zero when $m=0$.
Now we use  the connection (\ref{3.12}) between the Gauss hypergeometric function and the Chebyshev polynomials of the first kind  to obtain
\beq
J_{j}^{(\Ga)}(\Gz)=\cot\pi\Ga (1-\Gz)^{\Ga-1}(1+\Gz)^{-\Ga}T_j(\Gz)-q_{j-1}^{(\Ga)}(\Gz), \quad -1<\Gz<1, \quad 
j=0,1,\ldots,
\label{3.25}
\eeq
where $q_{-1}^{(\Ga)}(\Gz)\equiv 0$.
To finalize our computations, we substitute this expression into the first formula in (\ref{3.17}) to have
$I_j^{(\Ga)}(\Gz)$. According to (\ref{3.16}) the function $\Gf_j(x)$ is  the sum
of the functions
$$
\fr12(1+\Gz)^{\Gr_1}(1-\Gz)^{1-\Gr_1}I_j^{(\Gr_1)}(\Gz)=B(\Gz)+\fr12(1+\Gz)^{\Gr_1}(1-\Gz)^{1-\Gr_1}q_j^{(\Gr_1)}(\Gz),
$$
\beq
\fr12(1+\Gz)^{1-\Gr_1}(1-\Gz)^{\Gr_1}I_j^{(1-\Gr_1)}(\Gz)=-B(\Gz)+\fr12(1+\Gz)^{1-\Gr_1}(1-\Gz)^{\Gr_1}
q_j^{(1-\Gr_1)}(\Gz),
\label{3.25'}
\eeq
where
\beq
B(\Gz)=\fr12\cot\pi\Gr_1[C_j-T_{j+1}(\Gz)].
\label{3.25''}
\eeq
After some obvious simplifications this eventually yields
the spectral relation for the operator $\CS$ 
\beq
\CS[\Gf_j](x)=N_{j}-\cos(j+1)\pi x, \quad 0<x<1, \quad j=0,1,\ldots,
\label{3.26}
\eeq
where 
\beq
N_{2m+1}=0, \quad N_{2m}=\sum_{j=0}^{2m}\fr{(-2m)_j(2m)_j(\Gr_1)_j}{(1/2)_j(j!)^2},
 \quad m=0,1,\ldots,
 \label{3.27}
\eeq
and $\Gf_j(x)$ are expressed through the degree-$2j$ polynomials on $\sin\fr{\pi x}{2}$,
\beq
\Gf_j(x)=\cos^{2\Gr_1}\fr{\pi x}{2}\sin^{2(1-\Gr_1)}\fr{\pi x}{2} q_j^{(\Gr_1)}(x)+
\cos^{2(1-\Gr_1)}\fr{\pi x}{2}\sin^{2\Gr_1}\fr{\pi x}{2} q_j^{(1-\Gr_1)}(x),
 \label{3.27.0}
 \eeq
 where
 \beq
 q_j^{(\Ga)}(x)=\fr{1}{2\sin\pi\Ga}\sum_{m=1}^{j+1}\fr{(-j-1)_m(j+1)_m}{(1/2)_m m!}
 \sum_{k=0}^{m-1}\fr{(\Ga)_k}{k!}\sin^{2(m-1-k)}\fr{\pi x}{2}.
 \label{3.27.1}
 \eeq
 Evidently, this expression coincides with  (\ref{3.6}).

Note that  on putting $j=-1$ in (\ref{3.16})
and (\ref{3.25'}) we obtain
the relation
\beq
\fr{\sin\pi x}{2}\int_0^1\fr{[a^{2\Gr_1-1}(x,\tau)+a^{-2\Gr_1+1}(x,\tau)]d\tau}
 {\cos\pi\tau-\cos\pi x}=0, \quad 0<x<1,
 \label{3.27'}
 \eeq
 that is 
 \beq
 \CS^{-1}[1]=0,\quad  0<x<1.
\label{3.27''}
\eeq

\begin{figure}[t]
\centerline{
\scalebox{0.7}{\includegraphics{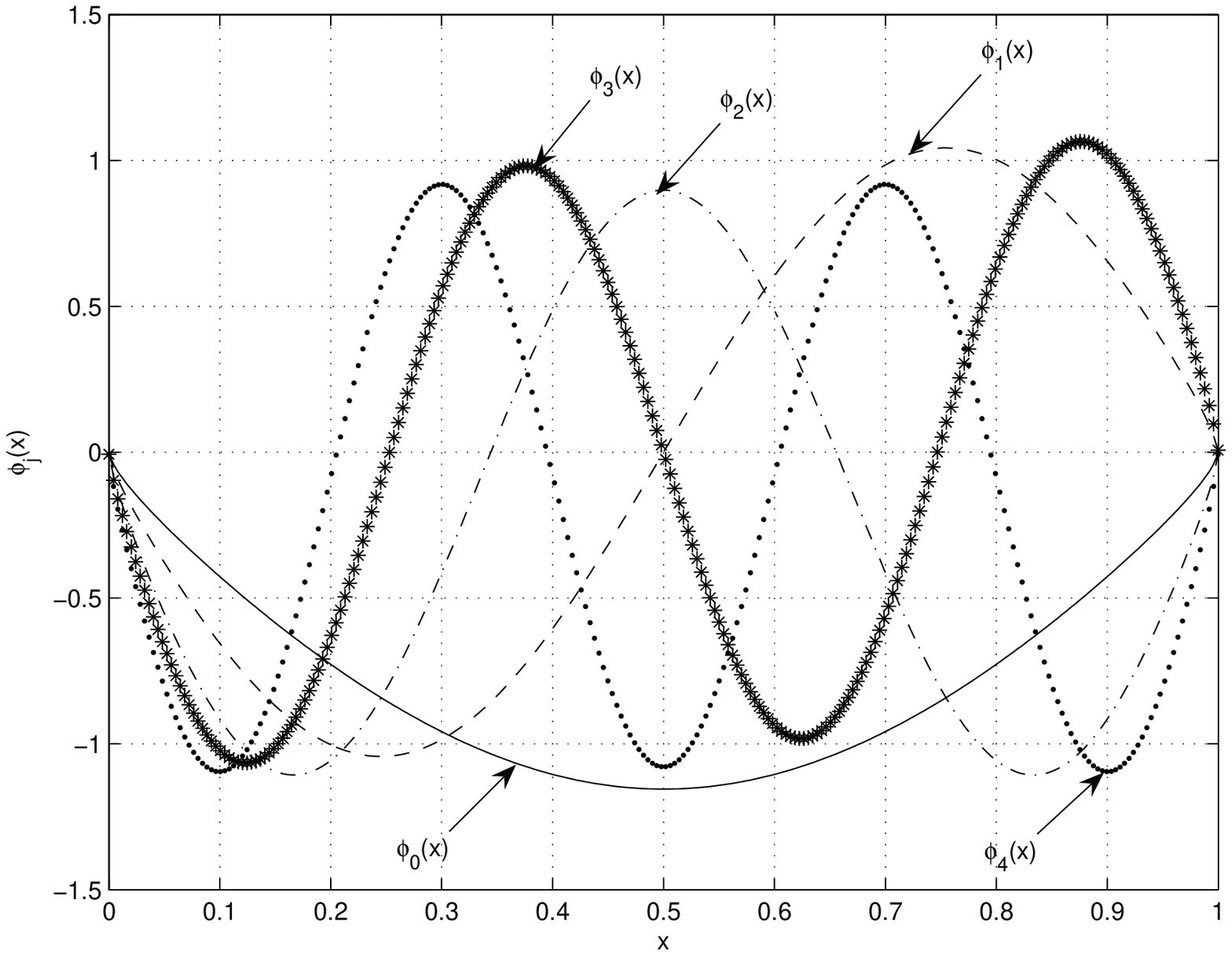}}}
\caption{Functions $\Gf_j(x)$, $j=0,1,\ldots,4$, $\Gb=0.5$}
\label{fig1}
\end{figure}  

In Figure \ref{fig1}, we plot the  functions $\Gf_j(x)$, $j=0,1,\ldots, 4$. 
They vanish at the endpoints $x=0$ and $x=1$. This numerical result is consistent with  formulas (\ref{2.37.3}).
Also, the functions $\Gf_j(x)$ have some properties which make them resemble 
orthogonal polynomials.  
Firstly,  the functions $\Gf_j(x)$ have exactly $j$ real roots
on the interval $(0,1)$. 
Secondly, although $\{\Gf_j(x)\}$, $j=0,1,\ldots,$ is not an orthogonal system, 
the  functions $\{\Gf_{2m+1}(x)\}$ are
 orthogonal  
to  $\{\Gf_{2n}(x)\}$ in the associated weight Hilbert space $L^2_w(0,1)$, 
 $w(x)=\sin^{2(\Gr_1-1)} \pi x$,
\beq
\int_0^1 \Gf_{2m+1}(x) \Gf_{2k}(x) w(x)dx=0, \quad m,k=0,1,\ldots.
\label{W.1}
\eeq
Finally, the functions $\Gf_j(x)$  ($j=0,1,\ldots$) can be employed for constructing a series-form 
solution of the characteristic equation in the same fashion as the Chebyshev polynomials 
are used for the integral equation with the Cauchy kernel in the segment. 

\subsection{Series-form solution of the characteristic equation}

We consider the characteristic equation 
\beq
\CS[\Gf](x)=C-F(x), \quad 0<x<1,
\label{W.3}
\eeq
with the constant $C$ given by 
\beq
 C=\fr{\sin\pi \Gr_1}{2}\int_0^1\left(\tan^{2\Gr_1-1}\fr{\pi x}{2}+\cot^{2\Gr_1-1}\fr{\pi x}{2}\right)F(x)dx
 \label{W.4}
 \eeq
and the function $F(x)$ represented by its cosine Fourier series 
$$
F(x)=f_0+2\sum_{j=1}^\infty f_j \cos\pi jx, \quad 0<x<1,
$$
\beq
f_j=\int_0^1 F(x)\cos\pi jx dx.
\label{W.5}
\eeq
According to Theorem 2.2 and since  $\CS^{-1}[1]=0$, we derive
$\Gf(x)=-\CS^{-1}[F](x)$, $0<x<1$. 
On replacing $F(x)$ by its cosine series and changing the order of integration and summation 
 we have
\beq
\Gf(x)=-2\sum_{j=1}^\infty f_j \CS^{-1}[\cos\pi j\Gx](x).
\label{W.6}
\eeq
We recall that $\CS^{-1}[\cos\pi j \Gx](x)=-\Gf_{j-1}(x)$, $j=1,2,\ldots.$ This brings us to
\beq
\Gf(x)=2\sum_{j=0}^\infty f_{j+1}\Gf_j(x).
\label{W.7}
\eeq
We consider two examples, $F(x)=x$ and $F(x)=x^2$.
The cosine Fourier series of these functions have the form
$$
x=\fr12+\fr{2}{\pi^2}\sum_{j=1}^\infty\fr{(-1)^j-1}{j^2}\cos\pi j x, \quad 0<x<1,
$$
\beq
x^2=\fr13+\fr{4}{\pi^2}\sum_{j=1}^\infty\fr{(-1)^j}{j^2}\cos\pi j x, \quad 0<x<1.
\label{W.8}
\eeq
We substitute these series into the relation (\ref{W.4}), use formula (\ref{3.15})
and determine that if $F(x)=x$, then $C=\fr12$ for all $\Gr_1$ ($|\Gb|<1$).  In the case $F(x)=x^2$,
\beq
C=\fr13+\fr{1}{\pi^2}\sum_{m=1}^\infty\GP_m,\quad
\GP_m=\fr{1}{m^2}\sum_{j=0}^{2m}\fr{(-2m)_j(2m)_j(\Gr_1)_j}{(1/2)_j(j!)^2}.
\label{W.9}
\eeq
Our computations show that for $\Gb=\fr12$, $C(M)=\fr13+\pi^{-2}\sum_{m=1}^M \GP_m \approx \fr23$
for large $M$: $C(100)=0.665659$, $C(1000)=0.666565$, $C(2000)=0.666616$,
$C(3000)=0.666633$. As for the solution itself, it turns out that the error of approximation 
decreases as $m_0$ in
$\Gf(x)\approx 2\sum_{j=0}^{m_0} F_{j+1}\Gf_j(x)$
approaches $20$, and
the algoritm becomes unstable for $m_0\ge 25$ (Table 1). This is a typical feature even in the case
of the classical method of orthogonal polynomials caused by their oscilation that increases as $m_0$ grows.
In Figure \ref{fig2}, we plot the  function $\Gf(x)$ when $m_0=20$ for the cases $F(x)=x$ and $F(x)=x^2$.
\vspace{1mm}
\begin{table}
	\centering
\begin{tabular}{|c|c|c|c|c|c|c|}
\hline
     & $m_0=5$ & $m_0=10$ & $m_0=15$ & $m_0=20$ & $m_0=22$    & $m_0=25$\\
\hline
$x_0=0.5$ & 0.445026  & 0.439400 & 0.440180& 0.441492& 0.441439  & 0.434170\\
\hline
$x=0.25$ & 0.371760  & 0.372910 & 0.370147& 0.371442& 0.371252  & 0.363680\\
\hline
\end{tabular}
	\caption{The values of the function $\Gf(x)$ as $\Gb=0.5$, $x=0.5$ and $x=0.25$ for some values of $m_0$.}
	\label{tab1}
\end{table}

\begin{figure}[t]
\centerline{
\scalebox{0.7}{\includegraphics{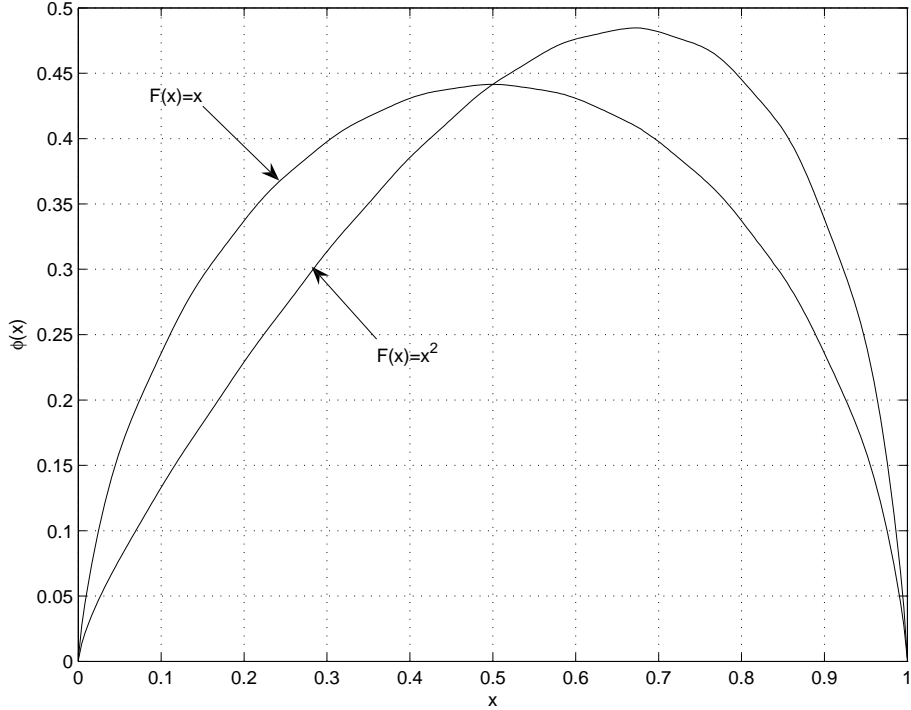}}}
\caption{The solution of the characteristic equation $\Gf(x)$ for $F(x)=x$ and $F(x)=x^2$ when  $\Gb=0.5.$}
\label{fig2}
\end{figure}

\subsection{Solution of the complete integral equation}

We continue considering the case $|\Gb|<1$.
In the class of functions bounded at the ends the complete singular integral equation (\ref{3.1})
is solvable if and only if the following condition is met:
\beq
\int_0^1\left(\tan^{2\Gr_1-1}\fr{\pi x}{2}+\cot^{2\Gr_1-1}\fr{\pi x}{2}\right)\left[C-F(x)-\int_0^1
K(x,\Gx)\Gf(\Gx)d\Gx\right]dx=0.
\label{3.28}
\eeq
Suppose  this condition holds. We build up the solution to the integral equation (\ref{3.1}) 
in the series form
\beq
\Gf(x)=\sum_{j=0}^\infty b_j\Gf_j(x),
\label{3.29}
\eeq
where the coefficients $b_j$ ($j=0,1,\ldots$) are to be determined. On substituting this series
into equation (\ref{3.1}), utilizing the expansion (\ref{3.29}) and the orthogonality relation 
\beq\int_0^1 \cos j\pi x \cos n\pi xdx=\left\{\begin{array}{cc}
1,& n=j=0,\\
\fr12, & n=j\ne 0,\\
0, & n\ne j,\\
\end{array}
\right.
\label{3.30}
 \eeq
we derive the following infinite system of algebraic equations:
\beq
\sum_{j=0}^\infty \left[N_{j+1}\Gd_{n,0}-\fr12\Gd_{n,j+1}+k_{nj}\right]b_j=
C\Gd_{n,0}-f_n,\quad n=0,1,\ldots,
\label{3.31}
\eeq
Here, $\Gd_{n, j}$ is the Kronecker delta, and
$$
 k_{nj}=\int_0^1\int_0^1 K(x,\Gx)\Gf_j(\Gx)\cos n\pi x d\Gx dx,
$$
\beq
f_n=\int_0^1 F(x)\cos n\pi x dx.
\label{3.32}
\eeq
It is convenient to split this system into 
two parts
\beq
n=0: \quad  \sum_{j=0}^\infty \left(N_{j+1}+k_{0j}\right)b_j=-f_0+C, 
\label{3.33}
\eeq
and 
 \beq 
n=1,2,\ldots: \quad -\fr12 b_{n-1}+\sum_{j=0}^\infty k_{nj}b_j=-f_n.
 \label{3.34}
 \eeq
It turns out that  equation (\ref{3.33}) is equivalent to the solvability relation (\ref{3.28}).
To prove this, we substitute the Fourier expansions
$$
\int_0^1 K(x,\Gx)\Gf_j(\Gx)d\Gx=k_{0j}+2\sum_{l=1}^\infty k_{nj}\cos \pi lx,
$$
\beq
F(x)=f_0+2\sum_{j=0}^\infty f_j\cos\pi jx,
\label{3.35}
\eeq
into equation (\ref{3.28}), recall formula (\ref{3.10}) and rewrite the solvability 
condition (\ref{3.28}) as
\beq
M_0\left(C-f_0-\sum_{j=0}^\infty k_{0j} b_j\right)-2\sum_{n=1}^\infty M_n\left(f_n+\sum_{j=0}^\infty 
k_{nj} b_j\right)=0.
\label{3.36}
\eeq
Now, the coefficients $b_n$ solve the system of equations (\ref{3.33}) and (\ref{3.34}).
This immediately brings us to the identity
\beq
\sum_{j=0}^\infty b_j M_{j+1}-\sum_{n=1}^\infty b_{n-1} M_n=0.
\label{3.37}
\eeq
and therefore,  in the class of functions bounded at the ends,
the singular integral equation (\ref{3.1}) is equivalent to the system
of linear algebraic equations (\ref{3.33}) and (\ref{3.34}); its solution automatically
satisfies the solvability condition (\ref{3.28}). In Appendix D, we use 
the classical method of orthogonal polynomials to analyze the complete singular integral equation with the Cauchy kernel in the class of bounded
at the endpoints functions and show that the corresponding solvability 
condition is also equivalent to the first ($n=0$) equation of the associated 
infinite system of algebraic equations.

\setcounter{equation}{0}

\section{Antiplane problem for a crack in a composite plane}\label{s4}

The problem under consideration is one of antiplane strain on a crack 
$0<x<1$, $y=0^\pm$.
The elastic medium is formed by two half-planes $x<0$ and $x>1$
and an infinite strip $0<x<1$. The shear moduli of the half-planes and
the strip are $G_1$ and $G_2$, respectively. The faces of the crack are
subjected to traction $\tau_{yz}=-f(x)$, $0<x<1$, $y=0^\pm$. The problem is governed by
the following boundary value problem for the Laplace operator in the plane:
$$
\GD w(x,y)=0, \quad |x|<\infty,\quad  |y|<\infty,\quad  x\ne 0, \quad  x\ne 1,
$$
$$
w(0^-,y)-w(0^+,y)=0, \quad w(1^-,y)-w(1^+,y)=0, \quad  |y|<\infty,
$$
$$
G_1w_x(0^-,y)=G_2w_x(0^+,y), \quad G_2w_x(1^-,y)=G_1w_x(1^+,y), \quad  |y|<\infty,
$$
\beq
G_2w_y(x,0^\pm)=-f(x),   \; 0<x<1.
\label{4.1}
\eeq
We note that due to the symmetry of the problem,
$w(x,0^\pm)=0$, $-\infty<x<0$, $1<x<\infty$, and the problem can be restated for say, the upper half-plane. On the crack faces,
the displacement $w$ is discontinuous, $w(x,0^+)=-w(x,0^-)$,
 $0<x<1$, and the displacement jump is to be determined. Denote
 $\Gf(x)=G_2w(x,0^+)$, $0<x<1$. By the method of integral transformations
 the problem reduces to the following integral equation {\bf(\ref{moi})}:
\beq
-\fr{d}{dx}\int_0^1 M(x,\Gx)\Gf(\Gx)d\Gx=f(x), \quad 0<x<1,
\label{4.2}
 \eeq
 where
 $$
 M(x,\Gx)=\fr{1}{\pi}\left[\fr{1}{\Gx-x}+\fr{\Gb}{\Gx+x}+\fr{\Gb}{x+\Gx-2}+R(x,\Gx)
\right],
$$$$
R(x,\Gx)=\Gb[D(x+\Gx)-D(2-x-\Gx)]+\Gb^2\left[D(2-x+\Gx)-D(2+x-\Gx)
+\fr{2(x-\Gx)}{4-(x-\Gx)^2}\right],
$$
\beq
D(x)=\sum_{j=1}^\infty\fr{\Gb^{2j}}{x+2j}, \quad \Gb=\fr{\Gl-1}{\Gl+1}\in(-1,1), \quad \Gl=\fr{G_1}{G_2}\in(0,\infty).
\label{4.3}
\eeq
By integrating equation (\ref{4.3}) with respect to $x$ we can rewrite
the new equation in the form used in the previous section
\beq
\int_0^1[S(x,\Gx)+K(x,\Gx)]\Gf(\Gx)d\Gx=-F(x)+C, \quad 0<x<1,
\label{4.4}
\eeq
where $F(x)=\int f(x)dx$ and
\beq
K(x,\Gx)= \fr{1}{\pi}\left[\fr{1}{\Gx-x}+\fr{\Gb}{\Gx+x}+\fr{\Gb}{x+\Gx-2}+R(x,\Gx)
\right]-\fr12\cot\fr{\pi(\Gx-x)}{2}-\fr{\Gb}{2}\cot\fr{\pi(\Gx+x)}{2}.
\label{4.6}
\eeq
is a regular kernel.
This equation has been solved in the previous section by reducing it to the infinite
system of linear algebraic equations of the second kind (\ref{3.34}). 
The coefficients (\ref{3.32}) of the infinite system can be represented in the 
form
$$
f_n=\fr{1}{\pi}\int_{-1}^1F\left(\fr{\cos^{-1}\Gz}{\pi}\right)\fr{T_n(\Gz)d\Gz}{\sqrt{1-\Gz^2}},
$$
\beq
k_{nj}=\fr{1}{\pi^2}\int_{-1}^1\int_{-1}^1
K\left(\fr{\cos^{-1}\Gz}{\pi}, \fr{\cos^{-1}\Gn}{\pi}\right)\Gf_j\left(\fr{\cos^{-1}\Gn}{\pi}\right)
\fr{T_n(\Gz)d\Gz d\Gn}{\sqrt{1-\Gz^2}\sqrt{1-\Gn^2}},
\label{4.7}
\eeq
and computed by the Gauss quadrature formulas
$$
f_n=\fr{1}{t_1}\sum_{m=1}^{t_1}F\left(\fr{2m-1}{2t_1}\right)\cos\fr{(2m-1)n\pi}{2t_1},
$$
\beq
k_{nj}=\fr{1}{t_1t_2}\sum_{m=1}^{t_1}\sum_{l=1}^{t_2}K\left(\fr{2m-1}{2t_1},\fr{2l-1}{2t_2}
\right)\Gf_j\left(\fr{2l-1}{2t_2}\right)\cos\fr{(2m-1)n\pi}{2t_1},
\label{4.8}
\eeq
where $t_1$ and $t_2$ are the numbers of the Gauss nodes. An approximate solution 
\beq
\Gf^{(N)}(x)=\sum_{j=0}^{N-1} b_j^{(N)}\Gf_j(x),
\label{4.9}
\eeq
 of the system (\ref{3.34})
 is found by the truncation method; the coefficients $b_j^{(N)}$
solve the system
\beq
-\fr12 b^{(N)}_{n-1}+\sum_{j=0}^{N-1} k_{nj}b_j^{(N)}=-f_n, \quad n=1,2,\ldots,N.
\label{4.10}
\eeq

 \begin{figure}[t]
\centerline{
\scalebox{0.6}{\includegraphics{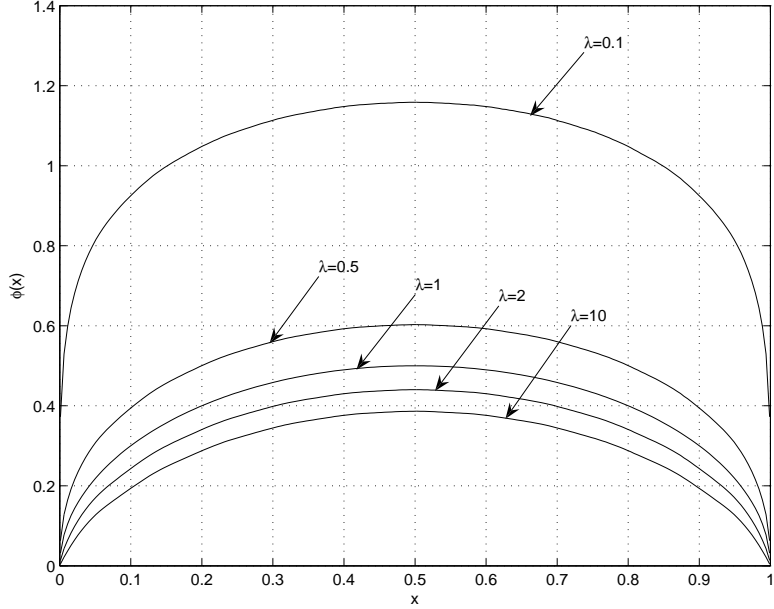}}}
\caption{Antiplane strain: the displacement function $\Gf(x)$ 
for some values of the parameter  $\Gl=G_1/G_2$.}
\label{fig3}
\end{figure} 

 \begin{figure}[t]
\centerline{
\scalebox{0.6}{\includegraphics{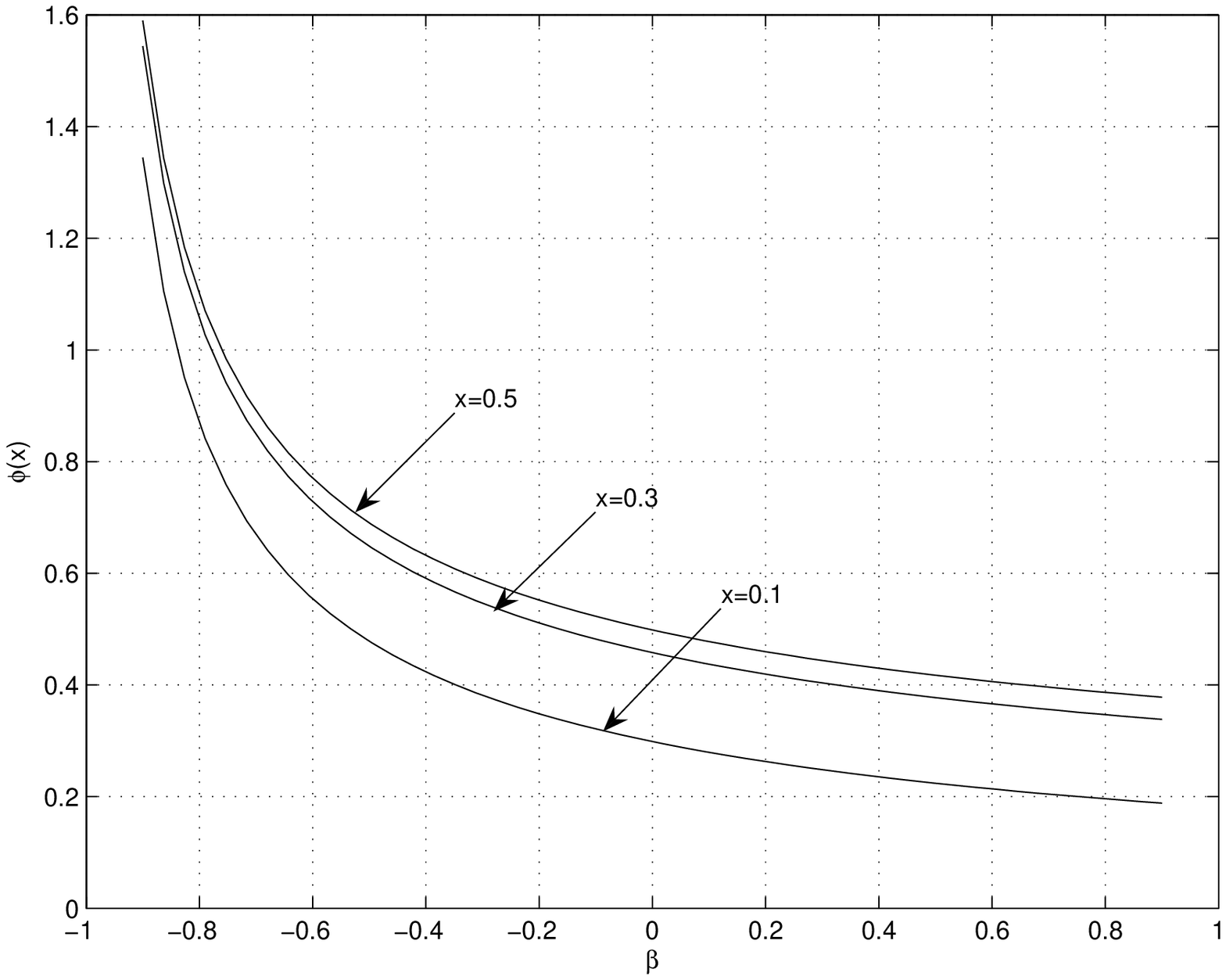}}}
\caption{Antiplane strain: the displacement function $\Gf(x)$ as a function of $\Gb=(G_1-G_2)/(G_1+G_2)$
for some values of $x$.}
\label{fig4}
\end{figure}

To test the efficiency of the numerical scheme, we consider the case of a uniform load, $f(x)=
P=\const$ or, equivalently, the case $F(x) = Px$ (for computations, we select $P=1$). 
Table 2 shows how the approximate solution $\Gf^{(N)}(x)$ given by (\ref{4.9}) 
depends on the numbers  $t_1$ and $t_2$ of the Gauss nodes when $x=0.5$, $\Gl=0.5$, $F(x)=x$, and $N=17$. 
In Table 3 we report the values $\Gf^{(N)}(0.5)$ when $\Gl=0.5$, $F(x)=x$, $t_1=200$, and $t_2=210$
for some values of the truncation parameter $N$. Again, as in the case of the characteristic equation,
when $N$ becomes greater than 20, then there is no improvement of the 
the accuracy of the numerical results, and  the error of approximations grows when $N\ge 25$. 
\vspace{1mm}
\begin{table}
	\centering
\begin{tabular}{|c|c|c|c|c|c|}
\hline
     & $t_1=10$ &   $t_1=50$ &    & $t_1=100$     & $t_1=300$\\
\hline
$t_2=11$ &    0.599598 & 0.0.600131&  $t_2=110$  & 0.601067 & 0.601099 \\
\hline
$t_2=51$ & 0.600073 & 0.600895 &$t_2=310$ &   0.601091 & 0.601123    \\
\hline

\end{tabular}
	\caption{Antiplane strain: the values of the function $\Gf^{(N)}(x)$ at $x=0.5$ for some numbers $t_1$ and $t_2$ of the Gauss nodes.}
	\label{tab2}
\end{table}

\vspace{1mm}
\begin{table}
	\centering
\begin{tabular}{|c|c|c|c|c|c|}
\hline
     & $N=5$ &   $N=10$ & $N=15$     & $N=20$  & $N=25$\\
\hline
$\Gf^{(N)}(0.5)$ &    0.582829 & 0.601814  & 0.602258 &  0.601812     & 0.604167 \\
\hline

\end{tabular}
	\caption{Antiplane strain: the values of the function $\Gf^{(N)}(x)$ at $x=0.5$ for some values of $N$.}
	\label{tab3}
\end{table}

Figure \ref{fig3} presents the computations for the
displacement of the points on the upper crack face, the function $\Gf(x)$, for some values of the parameter
$\Gl=G_1/G_2$. It is seen that the crack opening is growing when the parameter $\Gl$ is decreasing.
Figure \ref{fig4} shows that when $\Gb=(\Gl-1)(\Gl+1)^{-1}\to -1$, that is when $\Gl\to 0$, the crack opening is growing to infinity.
If $\Gb\to 1$, then $\Gl\to\infty$, and the crack opening tends to a certain limit, a function $\Gf_\infty(x)$.

 \setcounter{equation}{0}

\section{Plane strain of a composite plane with a crack}\label{s5}

 In this section, we generalize the algorithm to the biharmonic case.
As before, we consider a composite plane with a crack $0<x<1$, $y=0^\pm$.
The shear moduli and Poisson ratios of the half-planes $x<0$ and $x>1$
are the same, $G_1$ and $\nu_1$, respectively, while the corresponding elastic constants for the strip
$0<x<1$ are $G_2$ and $\nu_2$. The faces of the crack are subjected to loading
$\Gs_y=-f(x)$, $\tau_{xy}=0$, $0<x<1$, $y=0^\pm$.
The normal displacement $v$ is discontinuous across the crack faces, while the tangential
displacement $u$ is continuous. Denote
\beq
\Gf(x)=v(x,0^+)-v(x,0^-), \quad \supp \Gf(x)\subset [0,1].
\label{5.1} 
\eeq
Let $U(x,y)$ be the Airy function of the problem. Consider the conditions of plane strain.
Then  the function $U$ solves the following
discontinuous boundary value problem for the biharmonic operator:
$$
\GD^2 U(x,y)=0, \quad |x|<\infty, \quad |y|<\infty, \quad x\ne0, \quad x\ne 1,
$$$$
\fr{\Md^j U}{\Md y^j}(x, 0^+)- \fr{\Md^j U}{\Md y^j}(x, 0^-)=0, \; j=0,1,2,\quad
\fr{\Md^3 U}{\Md y^3}(x, 0^+)- \fr{\Md^3 U}{\Md y^3}(x, 0^-)=-\fr{2G_2}{1-\nu_2}\Gf''(x),\quad |x|<\infty,
$$$$
\fr{\Md^j U}{\Md x^j}(0^+,y)- \fr{\Md^j U}{\Md x^j}(0^-,y)=0, \quad 
\fr{\Md^j U}{\Md x^j}(1^+,y)- \fr{\Md^j U}{\Md x^j}(1^-,y)=0, \quad j=0,1, \quad |y|<\infty.
$$$$
\fr{1-\nu_1}{G_1}\fr{\Md^2 U}{\Md x^2}(0^-,y)
-\fr{1-\nu_2}{G_2}\fr{\Md^2 U}{\Md x^2}(0^+,y)=\left(\fr{\nu_1}{G_1}-\fr{\nu_2}{G_2}\right)
\fr{\Md^2 U}{\Md y^2}(0,y),\quad |y|<\infty.
$$
$$
\fr{1-\nu_1}{G_1}\fr{\Md^3 U}{\Md x^3}(0^-,y)
-\fr{1-\nu_2}{G_2}\fr{\Md^3 U}{\Md x^3}(0^+,y)=\left(\fr{\nu_1-2}{G_1}-\fr{\nu_2-2}{G_2}\right)
\fr{\Md^3 U}{\Md x\Md y^2}(0,y),\quad |y|<\infty,
$$$$
\fr{1-\nu_2}{G_2}\fr{\Md^2 U}{\Md x^2}(1^-,y)
-\fr{1-\nu_1}{G_1}\fr{\Md^2 U}{\Md x^2}(1^+,y)=\left(\fr{\nu_2}{G_2}-\fr{\nu_1}{G_1}\right)
\fr{\Md^2 U}{\Md y^2}(1,y),\quad |y|<\infty,
$$$$
\fr{1-\nu_2}{G_2}\fr{\Md^3 U}{\Md x^3}(1^-,y)
-\fr{1-\nu_1}{G_1}\fr{\Md^3 U}{\Md x^3}(1^+,y)=\left(\fr{\nu_2-2}{G_2}-\fr{\nu_1-2}{G_1}\right)
\fr{\Md^3 U}{\Md x\Md y^2}(1,y),\quad |y|<\infty,
$$
\beq
\fr{\Md^2 U}{\Md x^2}(x,0^\pm)=-f(x), \quad 0<x<1.
\label{5.2}
\eeq
It is natural to apply the Fourier transform with respect to $y$ 
\beq 
U_\Ga(x)=\int_{-\infty}^\infty U(x,y) e^{i\Ga y}dy
\label{5.3}
\eeq
and reduce the problem
to the one-dimensional discontinuous boundary value problem 
$$
\left(\fr{d^4}{dx^4}-2\Ga^2\fr{d^2}{dx^2}+\Ga^4\right)U_\Ga(x)=-\fr{2G_2}{1-\nu_2}\Gf''(x), \quad |x|<\infty,\quad x\ne 0, \quad x\ne 1,
$$$$
\fr{d^j}{dx^j}U_\Ga(0^-)=\fr{d^j}{dx^j}U_\Ga(0^+), \quad \fr{d^j}{dx^j}U_\Ga(1^-)=\fr{d^j}{dx^j}U_\Ga(1^+),
$$$$
\fr{1-\nu_1}{G_1}\fr{d^{j+2}}{dx^{j+2}}U_\Ga(0^-)-\fr{1-\nu_2}{G_2}\fr{d^{j+2}}{dx^{j+2}}
U_\Ga(0^+)+\Ga^2\Gl_j \fr{d^{j}}{dx^{j}}U_\Ga(0)=0,
$$
\beq
\fr{1-\nu_2}{G_2}\fr{d^{j+2}}{dx^{j+2}}U_\Ga(1^-)-\fr{1-\nu_1}{G_1}\fr{d^{j+2}}{dx^{j+2}}
U_\Ga(1^+)-\Ga^2\Gl_j \fr{d^{j}}{dx^{j}}U_\Ga^{j}(1)=0,
\quad j=0,1,
\label{5.4}
\eeq
where
\beq
\Gl_0=\fr{\nu_1}{G_1}-\fr{\nu_2}{G_2}, \quad \Gl_1=\fr{\nu_1-2}{G_1}-\fr{\nu_2-2}{G_2}.
\label{5.5}
\eeq
Associated with the differential operator in (\ref{5.4}) is the fundamental function 
\beq
\GF_\Ga(x,\Gx)=\fr{1+|\Ga||x-\Gx|}{4|\Ga|^3}e^{-|\Ga||x-\Gx|}.
\label{5.6}
\eeq
After integration by parts two times the general solution of the problem (\ref{5.4})  becomes
$$
U_\Ga(x)=\fr{G_2}{2(1-\nu_2)|\Ga|}\int_0^1(1-|\Ga||x-\Gx|)  e^{-|\Ga||x-\Gx|}\Gf(\Gx)d\Gx
$$
\beq
+\left\{\begin{array}{cc}
(c_{00}+c_{01}x) e^{|\Ga|x}, & x<0,\\
(c_{10}+c_{11}x) \cosh |\Ga|x+(c_{12}+c_{13}x) \sinh |\Ga|x, & 0<x<1,\\
(c_{20}+c_{21}x) e^{-|\Ga|x}, & x>1.\\
\end{array}
\right.
\label{5.7}
\eeq
The eight arbitrary constants in this solution are fixed by the eight conditions 
at the points $x=0$ and $x=1$ in (\ref{5.4}). By inversion of the Fourier transform and 
satisfying the last condition in (\ref{5.2}) equivalent to
\beq
\fr{d}{dx}\fr{1}{2\pi}\int_{-\infty}^{\infty}U_\Ga'(x)e^{-i\Ga y} d\Ga|_{y=0}=-f(x), \quad 0<x<1,
\label{5.8}
\eeq
we eventually arrive at the following singular integral equation with two fixed singularities at the ends:
$$
\fr{1}{\pi}\int_0^1\left[
\fr{1}{\Gx-x}+\fr{b_1\Gx^2+b_2\Gx x+b_3 x^2}{(\Gx+x)^3}
+\fr{b_1(\Gx-1)^2+b_2(\Gx-1) (x-1)+b_3 (x-1)^2}{(\Gx+x-2)^3}
\right.
$$
\beq
\left.
+K_0(x,\Gx)\right]\Gf(\Gx)d\Gx=
-F(x)+C, \quad 0<x<1.
\label{5.9}
\eeq
Here,
$$
b_1=\fr{1}{\Gd_0}[(\nu_0+\mu_0-1)^2-4(1-\mu_0^2)],
\quad
b_2=\fr{4}{\Gd_0}[\nu_0(\nu_0-2)-3(1-\mu_0^2)],
$$$$
b_3=\fr{1}{\Gd_0}[-4\nu_0(\nu_0-2)+3(\nu_0+\mu_0-1)^2],\quad 
\Gd_0=(3+\mu_0-\nu_0)(1+3\mu_0+\nu_0),
$$
\beq
\mu_0=\fr{G_1(1-\nu_2)}{G_2(1-\nu_1)},\quad
\nu_0=\fr{\nu_1}{1-\nu_1}-\mu_0\fr{\nu_2}{1-\nu_2},
\label{5.10}
\eeq
$C$ is an arbitrary constant due to integration of equation (\ref{5.8}),
\beq
F(x)=\fr{4(1-\nu_2)}{G_2}\int f(x)dx,
\label{5.11}
\eeq
and the function $K_0(x,\Gx)$ is a regular kernel whose representation is quite
complicated and omitted.

Since the structure of the singular kernel in equation (\ref{5.9}) is different from the one 
in the governing  equation (\ref{4.4}) for the antiplane problem, we cannot directly
apply the algorithm of  section \ref{s3}. To adjust the scheme  to the plane case, 
we build up a new singular operator that generates a solution with the same   
singularities at the endpoints.
Assume that $\Gf(x)\sim A x^\Gg$, $x\to 0^+$, $A$ is a nonzero constant.
Because the function $\Gf(x)$ is the displacement jump of the normal displacement
across the crack, $\R\Gg$ is positive. To deal with the largest class of functions possible, we assume
 $\R\Gg\in(0,1)$. 
Then we employ the formulas
$$
\fr{1}{\pi}\int_0^1\fr{\Gx^\Gg d\Gx}{\Gx-x}= -\cot\pi\Gg x^\Gg+\Go_0(x), \quad x\to 0^+,
$$
\beq
\fr{1}{\pi}\int_0^1\fr{\Gx^{\Gg+2} d\Gx}{(\Gx+x)^3}=  -\fr{(\Gg+1)(\Gg+2) x^\Gg}{2\sin\pi\Gg}+\Go_1(x), \quad x\to 0^+,
\label{5.12}
\eeq
where $\Go_j(x)$ are continuously differentiable functions in the segment $0\le x<\Gve$
for some positive $\Gve$, and $\Go_j(0)\ne 0$, $j=0,1$.
Analysis of the singular integrals in (\ref{5.9})  as $x\to 0^+$ leads to the following
transcendental equation for the parameter $\Gg$:
\beq 
\GL(\Gg)\equiv \Gd_0\cos\pi\Gg-2[\mu_0^2-3-2\mu_0(\nu_0-1)+\nu_0(\nu_0-2)]\Gg^2-4(1-\mu_0^2)+(\nu_0+\mu_0-1)^2=0.
\label{5.13}
\eeq
\begin{figure}[t]
\centerline{
\scalebox{0.6}{\includegraphics{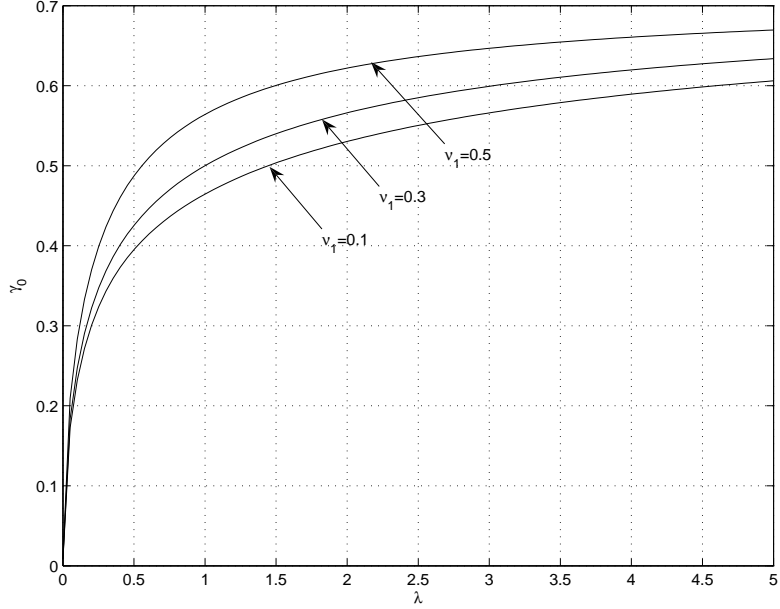}}}
\caption{The zero $\Gg_0\in(0,1)$ of the equation (\ref{5.13}) as a function of $\Gl=G_1/G_2$ 
for some values of $\nu_1$ when $\nu_2=0.3$.}
\label{fig5}
\end{figure} 
It turns out that in the strip $0<\R\Gg<1$ the function $\GL(\Gg)$ has one and only one zero, $\Gg_0$,
and it is real.
Its dependence on the parameter $\Gl=G_1/G_2$   
for some values of the Poisson ratio $\nu_1$ when $\nu_2=0.3$  is shown in Figure \ref{fig5}. It is seen that
when $\nu_1=\nu_2=0.3$ and $\Gl=1$, that is when the plane is homogeneous,  the parameter $\Gg_0$ is equal to $1/2$.
Also, if the shear modulus of the internal strip is greater 
than that of the surrounding matrix ($\Gl<1$) and  $\nu_1=\nu_2$, then $\Gg_0<1/2$, and
$\Gg_0\to 0$ as $\Gl\to 0$.   
On the other hand, if $\Gl\to\infty$, then the parameter
$\Gg_0\to\Gg_\infty$, where $\Gg_\infty$ is independent of $\nu_1$ and $\nu_2$ and $\Gg_\infty\approx 0.7111773$.

Introduce now a new singular integral equation 
\beq
\fr{1}{\pi}\int_0^1\left[\fr{1}{\Gx-x}+\fr{\Gb}{\Gx+x}+\fr{\Gb}{x+\Gx-2}\right]\psi(\Gx)d\Gx=-F(x)+C, \quad 0<x<1,\quad \Gb=-\cos\pi\Gg_0,
\label{5.14}
\eeq
associated with the complete singular equation (\ref{5.9}). Analysis of the singular integrals
in (\ref{5.14}) shows that
the derivatives of the functions $\Gf(x)$ and
$\psi(x)$ have identical asymptotic representations at the ends. We have
$$
\Gf'(x)\sim\psi'(x)\sim C_0 x^{\Gg-1}, \quad x\to 0^+,
$$
\beq
\Gf'(x)\sim\psi'(x)\sim C_1 (1-x)^{\Gg-1}, \quad x\to 1^-,
\label{5.15}
\eeq
where $C_j\ne 0$, $j=0,1$.
Therefore it is natural to rewrite the integral equation (\ref{5.9}) in the form
\beq
\fr{1}{\pi}\int_0^1\left[\fr{1}{\Gx-x}+\fr{\Gb}{\Gx+x}+\fr{\Gb}{x+\Gx-2}+R(x,\Gx)\right]\Gf(\Gx)d\Gx=-F(x)+C, \quad 0<x<1,
\label{5.16}
\eeq
where 
$$
R(x,\Gx)=\fr{b_1\Gx^2+b_2\Gx x+b_3 x^2}{(\Gx+x)^3}
+\fr{b_1(\Gx-1)^2+b_2(\Gx-1) (x-1)+b_3 (x-1)^2}{(\Gx+x-2)^3}
$$
\beq
+\fr{\cos\pi\Gg_0}{\Gx+x}+\fr{\cos\pi\Gg_0}{x+\Gx-2}+K_0(x,\Gx),
\label{5.17}
\eeq
which can be recast as equation (\ref{4.4}) with the kernel (\ref{4.6}), the function $R(x,\Gx)$ given by (\ref{5.17}) and the parameter $\Gb=-\cos\pi\Gg_0$.
An approximate solution of this equation is constructed 
in section \ref{s3}.

\begin{figure}[t]
\centerline{
\scalebox{0.6}{\includegraphics{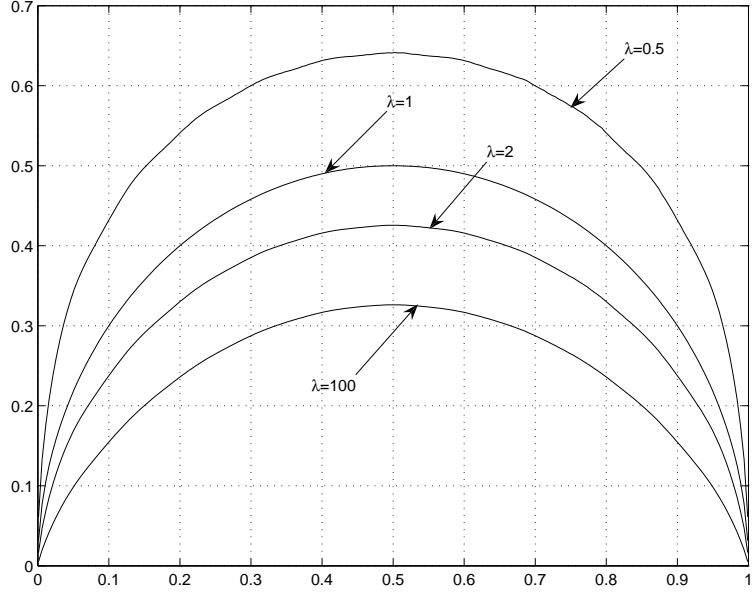}}}
\caption{Plane strain: the solution to the singular equation (\ref{5.18}), the function $\Gf_(x)$,
for some values of the parameter $\Gl=G_1/G_2$ when $F(x)=x$, $\nu_1=\nu_2=0.3$.}
\label{fig6}
\end{figure} 

To verify the technique, we drop the regular part $K_0(x,\Gx)$ and derive an approximate solution to the equation 
$$
\fr{1}{\pi}\int_0^1\left[
\fr{1}{\Gx-x}+\fr{b_1\Gx^2+b_2\Gx x+b_3 x^2}{(\Gx+x)^3}
+\fr{b_1(\Gx-1)^2+b_2(\Gx-1) (x-1)+b_3 (x-1)^2}{(\Gx+x-2)^3}\right]\Gf(\Gx)d\Gx=
$$
\beq
-F(x)+C, \quad 0<x<1,
\label{5.18}
\eeq
when $f(x)=P$ and $P=[4(1-\nu_2)]^{-1}G_2$. In this case $F(x)=x$. The function $\Gf(x)$ is plotted
in Figure \ref{fig6} for some values of the parameter $\Gl=G_1/G_2$ when $\nu_1=\nu_2=0.3$. The algorithm 
provides a good accuracy for all values of $\Gl$. Our computations show (Table 4) that for small values of the parameter $\Gl$
the same accuracy requires a larger dimension $N$ of the truncated system   
(\ref{4.10}).

\vspace{1mm}
\begin{table}
	\centering
\begin{tabular}{|c|c|c|c|c|c|}
\hline
     & $N=5$ &   $N=10$ & $N=15$     & $N=20$  & $N=25$\\
\hline
$\Gl=0.3$ &    0.744136 & 0.785573  & 0.787570 &  0.786676     & 0.782490 \\
\hline
$\Gl=0.5$ &    0.612515 & 0.640310  & 0.640704 &  0.639210     & 0.637523 \\
\hline
$\Gl=2$ &    0.414315 & 0.426399  & 0.425883 &  0.424791     & 0.424394 \\
\hline
$\Gl=100$ &    0.320104 & 0.326844  & 0.326405 &  0.325752    & 0.325526 \\
\hline

\end{tabular}
	\caption{Plane strain: the values of the function $\GF^{(N)}(x)$ at $x=0.5$ for some values of $N$ and $\Gl$
	when $\nu_1=\nu_2=0.3$.}
	\label{tab4}
\end{table}

\setcounter{equation}{0}

\section{Conclusions}

We have analyzed the singular integral equation $\CS[\Gf](x)=f(x)$, $0<x<1$, 
in the class of functions bounded at the ends. 
 The singular operator $\CS$ is defined by 
\beq
\CS[\Gf](x)=\int_0^1
\left[\fr12\cot\fr{\pi(\Gx-x)}{2}+\fr{\Gb}{2}\cot\fr{\pi(\Gx+x)}{2}
\right]\Gf(\Gx)d\Gx,
\label{6.1}
\eeq
and its kernel has fixed singularities at the points $x=\Gx=0$ and $x=\Gx=1$.
By reducing it to a vector Riemann-Hilbert problem with a piece-wise constant
matrix  coefficient we have found that, in general, in the class of H\"older functions
bounded at the ends, the solution does not exist. If a certain integral condition
is satisfied, then the solution exists, it is unique and given by a quadrature.
We have shown that in the case $|\Gb|<1$ the solution is monotonically decaying at the ends.
If $\Gb>1$, the solution oscillates at the ends and vanishes.
If $\Gb<-1$, then there are two possibilities to derive a solution.
Each one requires a certain condition of solvability and gives a solution that oscillates 
and does not vanish at one end and oscillates and vanishes at the second end.

We have managed to  obtain a spectral relation
$\CS[\Gf_j](x)=N_{j+1}-\cos[(j+1)\pi x]$, $0<x<1$, $j=0,1,\ldots$,
 for the singular operator $\CS$. Here, the numbers $N_j$ are given by (\ref{3.27}), the functions $\Gf_j(x)$
are
 \beq
\Gf_j(x)=\cos^{2\Gr_1}\fr{\pi x}{2}\sin^{2(1-\Gr_1)}\fr{\pi x}{2} q_j^{(\Gr_1)}(x)+
\cos^{2(1-\Gr_1)}\fr{\pi x}{2}\sin^{2\Gr_1}\fr{\pi x}{2} q_j^{(1-\Gr_1)}(x),
 \label{6.2}
 \eeq
 with $q_j^{(\Ga)}(x)$ being the degree-$j$ trigonometric polynomials
 \beq
 q_j^{(\Ga)}(x)=\sum_{\nu=0}^jc_{j\nu}^{(\Ga)}\left(\fr{1-\cos\pi jx}{2}\right)^\nu,\quad
 c_{j\nu}^{(\Ga)}=
 \fr{1}{2\sin\pi\Ga}\sum_{m=\nu+1}^{j+1}\fr{(-j-1)_m(j+1)_m(\Ga)_{m-1-\nu}}{(1/2)_m m!(m-1-\nu)!}.
 \label{6.3}
 \eeq
This  spectral relation has been used as the key step in the approximate scheme
for the complete singular integral equation with two fixed singularities. On expanding the 
unknown function in terms of the functions $\Gf_j(x)$ we have reduced the integral
equation to an infinite system of linear algebraic equations of the second kind
and solved  it by the reduction method.

The method has been applied to the antiplane problem for a finite crack in a composite
plane when the crack is orthogonal to the interfaces between a strip and two half-planes.
The crack lies in the strip, and its tips fall in the interfaces. The shear moduli
are the same for the half-planes, while the shear modulus of the strip is different. The problem
is governed by a complete singular integral equation with two fixed singularities,
and the generalized method of orthogonal polynomials proposed has been applied. The 
numerical algorithm has been successfully tested; it has a good  accuracy and it is rapidly convergent. We have further modified the method to adjust it to the solution
of the singular  integral equation with two fixed singularities arising in biharmonic problems.
We have derived the governing singular integral equation for the 
plane strain problem with the same geometry as in the antiplane case.
We have shown that if the singularities of the solution at the endpoints are real,
then it is possible to replace the singular operator associated with the plane problem 
by a simpler operator which satisfies the spectral relation used in the antiplane case.
A numerical test for the dominant singular integral equation associated with the plane strain
problem has been implemented. A good accuracy and fast convergence of the algorithm
has  been achieved.

\vspace{.2in}

{\centerline{\Large\bf  References}}

\begin{enumerate}

\item\label{klu}
P.I. Klubin, 
The calculation of girder and circular plates on an elastic foundation,
{\it Akad. Nauk SSSR. In\v{z}enerny\u{i} Sbornik} {\bf12} (1952) 95-135.

\item\label{pop1} G.Ia. Popov, Some properties of classical polynomials and their application to contact problems, {\it J. Appl. Math. Mech. } {\bf 27} (1963) 1255-1271.

\item\label{pop2} G.Ia. Popov,  On the method of orthogonal polynomials in contact problems of the theory of elasticity, {\it J. Appl. Math. Mech.} {\bf 33} (1969) 503-517.

\item\label{erd} F. Erdogan, Approximate solutions of systems of singular integral equations, 
{\it SIAM J. Appl. Math.} {\bf 17} (1969) 1041-1059.

\item\label{ant}
Y.A.Antipov, Nonlinear bending models for beams and plates, {\it Proc. R. Soc. A.} (2014) 470,
 20140064; DOI: 10.1098/rspa.2014.0064.

\item\label{moi}
N.G. Moiseyev and G.Ya. Popov, The antiplane problem of a crack with edges touching planes where the constants of elasticity change, {\it J. Appl. Math. Mech.}
{\bf 58} (1994) 713-725.

\item\label{gak} F.D. Gakhov, {\it Boundary Value Problems} (Pergamon Press, Oxford 1966).

\item\label{mus} N.I. Muskhelishvili, {\it Singular Integral Equations} (P. Noordhoff,
Groningen 1953). 

\item\label{bat2} H. Bateman, {\it Tables of Integral Transforms}, Vol.2 (McGraw-Hill Book Company, New York
1954).

\item\label{gra}  I. S. Gradshteyn and I. M. Ryzhik, {\it Table of Integrals, Series, and Products}
(Academic Press, Amsterdam 2007.

\item\label{bat} H. Bateman, {\it Higher Transcendental Functions}, Vol.1 (McGraw-Hill Book Company, New York
1953).

\end{enumerate}

\vspace{.1in}

\appendix

\setcounter{equation}{0}

\section{Proof of Theorem 2.1}

To prove the first part of Theorem 2.1, we 
follow  {\bf(\ref{moi})} and introduce a vector-function 
\beq
\BGY(z)=\left(
\begin{array}{c}
\Gvf_1(z)\\
\Gvf_2(z)\\
\end{array}
\right) 
\label{A.1}
\eeq
analytic in the semi-disc $D$ and H\"older continuous up to the boundary $L\cup (-1,1)$.
To extend its definition into the whole plane, we continue analytically this vector first into the lower half-disc by the relation
\beq
\BGY(z)=\left(
\begin{array}{cc}
0 & 1\\
1 & 0\\
\end{array}
\right)\ov{\Psi(\ov{z})}, \quad \I z < 0, \quad |z|<1, 
\label{A.2} 
\eeq
and then into the exterior of the disc $|z|<1$ by the law
\beq
\BGY(z)=-\ov{\BGY\left(\fr{1}{\ov{z}}\right)}, \quad |z|>1.
\label{A.3}
\eeq
Notice that then the vectors
\beq
\BGY^+(t)=\lim_{z\to t\in L, z\in D}\BGY(z), \quad  \BGY^-(t)=-\lim_{z\to t\in L, z\in D}\ov{\BGY(\ov{z^{-1}})},
\label{A.4}
\eeq
admit analytic continuation from the contour $L$ into the domains $D$ and
$\{|z| >1, \I z>0\}$, respectively.
Likewise, the vectors 
\beq
\BGY^+(t)=\lim_{z\to t\in(\infty, +\infty), z\in{\Bbb C}^+}\Psi(z), \quad  \BGY^-(t)=\lim_{z\to t\in(\infty, +\infty), z\in{\Bbb C}^+ }
\left(\begin{array}{cc}
0 & 1\\
1 & 0\\
\end{array}\right)\ov{\BGY(\ov{z})}
\label{A.5}
\eeq
admit analytic continuation from the real axis into the upper 
and lower half-planes, ${\Bbb C}^+$ and ${\Bbb C}^-$, respectively. We now invoke the boundary conditions 
(\ref{2.9}) to derive a Riemann-Hilbert problem for
the vector $\Psi(z)$ with a piecewise constant matrix coefficient.
Its boundary condition reads
\beq
\BGY^+(t)=A(t)\BGY^-(t)+\Bb(t), \quad t\in L\cup (-\infty,+\infty),
\label{A.6}
\eeq
where
$$
A(t)=\left\{\begin{array}{cc}
A_0, & t\in(-\infty,+\infty),\\
I, & t\in L,\\
\end{array}
\right.
\quad
A_0=\left(\begin{array}{cc}
1-\Gb & \Gb\\
-\Gb & 1+\Gb\\
\end{array}
\right),
\quad
I=\left(\begin{array}{cc}
1 & 0\\
0 & 1\\
\end{array}
\right),
$$
\beq
\Bb(t)=\left\{\begin{array}{cc}
{\bf 0}, & t\in(-\infty,+\infty),\\
\Bb_0(t), & t\in L,\\
\end{array}
\right.\quad 
{\bf 0}=\left(\begin{array}{c}
0\\
0 \\
\end{array}
\right),
\quad \Bb_0(t)=\left(\begin{array}{c}
0\\
2u(t) \\
\end{array}
\right).
\label{A.7}
\eeq
The matrix $A(t)$ can be factorized as
\beq
A(t)=X^+(t)[X^-(t)]^{-1}, \quad t\in L\cup(-\infty,+\infty),
\label{A.8}
\eeq
where
\beq
X(z)=\left(\begin{array}{cc}
1-\Gg(z) & \Gg(z)\\
-\Gg(z) & 1+\Gg(z)
\end{array}
\right),
\label{A.9}
\eeq
The function $\Gg(z)$ solves the problem
\beq
\Gg^+(t)-\Gg^-(t)=\Gb,\quad t\in L\cup(-\infty,+\infty),
\label{A.10}
\eeq
and has the form
\beq
\Gg(z)=\fr{\Gb}{2\pi i}[\log z-\log(-z)], \quad -\pi\le\arg z\le\pi.
\label{A.11}
\eeq
Because of the following properties of the function $\Gg(z)$: 
$$
\Gg^\pm(t)=\pm\fr{\Gb}{2}, \quad t\in(-\infty,+\infty),\qquad \Gg^+(t)=\Gg^-(t), \quad t\in L,
$$ 
\beq
\Gg(z)=-\ov{\Gg(\ov{z})}, \quad \I z<0, \qquad \Gg(z)=\ov{\Gg\left(\fr{1}{\ov{z}}\right)},\quad |z|>1,
\label{A.12}
\eeq
  the general solution of the vector Riemann-Hilbert problem (\ref{A.6}) in the class of symmetric functions (\ref{A.2}) and (\ref{A.3})
has the form
\beq
\Gvf_j(z)=\fr{1}{2\pi i}\int_L W_j(z,\tau)\fr{u(\tau)d\tau}{\tau}+ic[(-1)^j+2\Gg(z)], \quad j=1,2
\label{A.13}
\eeq
where $c$ is an arbitrary real constant and
$$
W_1(z,\tau)=
[\Gg(z)-\Gg(\tau)]\fr{\tau+z}{\tau-z}+[1-\Gg(z)-\Gg(\tau)]\fr{1+\tau z}{1-\tau z},
$$
\beq
W_2(z,\tau)=
[1+\Gg(z)-\Gg(\tau)]\fr{\tau+z}{\tau-z}-[\Gg(z)+\Gg(\tau)]\fr{1+\tau z}{1-\tau z}.
\label{A.14}
\eeq
It is directly verified that
$$
\ov{W_{3-j}(\ov{z},\tau)}=W_j(z,\tau), \quad \I z<0,
$$
\beq
\ov{W_j\left(\fr{1}{\ov{z}},\tau\right)}=-W_j(z,\tau),\quad |z|>1,\quad \tau\in L,\quad j=1,2,
\label{A.15}
\eeq
and therefore the conditions (\ref{A.2}) and (\ref{A.3}) are fulfilled.

Next we determine $\I\Gvf_2(t)$, $t\in L$, from (\ref{A.13})
\beq
\I\Gvf_2(t)=c[1+2\Gg(t)]-\fr{1}{2\pi}\int_L\left[\fr{\tau+t}{\tau-t}-\fr{\Gb(1+\tau t)}{1-\tau t}\right]\fr{u(\tau)d\tau}{\tau},\quad t\in L.
\label{A.16}
\eeq
Henceforth the function $u(\tau)=\R\Gvf_2(\tau)$ solves the integral equation (\ref{2.4})
if $c=0$ and $\I\Gvf_2(t)=-v(t)$. Finally notice that 
since
\beq
\lim_{z\to 0}[\Gvf_1(z)-\Gvf_2(z)]=-2ic,
\label{A.17}
\eeq
the constant $c$ vanishes
if  the functions in (\ref{A.13}) meet the condition (\ref{2.10}).

The inverse statement of Theorem 2.1 is proved on the basis of the Sokhotski-Plemelj formulas
and the representations (\ref{2.8}).

\setcounter{equation}{0}

\section{Cases $\Gb=\pm 1$}\label{spm1}

If $\Gb=\pm 1$, then we extend the definition of the functions $\Gf(x)$ and $f(x)$ to the interval 
$(-1,0]$ by the relations
\beq
\Gf(x)=-\Gb\Gf(-x),\quad f(x)=\Gb f(-x), \quad -1<x\le 0,
\label{2.78}
\eeq
 and convert equation (\ref{2.1}) into the equation
\beq
\fr12\int_{-1}^1\cot\fr{\pi(\Gx-x)}{2}\Gf(\Gx)d\Gx=f(x), \quad -1<x<1.
\label{2.79}
\eeq
By making the substitutions $\Gs=\pi(\Gx+1)$ and $s=\pi(x+1)$ and denoting
$\Gf(\Gx)=\hat\Gf(\Gs)$ and $f(x)=\hat f(s)$ we arrive at the equation
\beq
\fr{1}{2\pi}\int_0^{2\pi}\cot\fr{\Gs-s}{2}\hat\Gf(\Gs)d\Gs=\hat f(s),
\quad 0<s<2\pi.
\label{2.80}
\eeq
It is known {\bf(\ref{gak})}, pp.44, 244 that it is solvable if and only if
\beq
\int_0^{2\pi}\hat f(s)ds=0,
\label{2.81}
\eeq
and its solution is given by
\beq
\hat\Gf(s)=-\fr{1}{2\pi}\int_0^{2\pi}\hat f(\Gs)\cot\fr{\Gs-s}{2}d\Gs+C, \quad C=\fr{1}{2\pi}\int_0^{2\pi}\hat\Gf(\Gs)d\Gs.
\label{2.82}
\eeq
Therefore, in the case $\Gb=1$, the condition
\beq
\int_0^1 f(\Gx)d\Gx=0
\label{2.83}
\eeq
is necessary and sufficient for equation (\ref{2.1}) to be solvable.
If it is satisfied, the solution  is unique and has the form
\beq
\Gf(x)=-\fr12\int_0^1
\left[\cot\fr{\pi(\Gx-x)}{2}-\cot\fr{\pi(\Gx+x)}{2}
\right]f(\Gx)d\Gx, \quad 0<x<1.
\label{2.84}
\eeq
In the case $\Gb=-1$, the solution always exists, and its solution is defined up to an arbitrary constant
\beq
\Gf(x)=-\fr12\int_0^1
\left[\cot\fr{\pi(\Gx-x)}{2}+\cot\fr{\pi(\Gx+x)}{2}
\right]f(\Gx)d\Gx+C, \quad 0<x<1.
\label{2.85}
\eeq

\setcounter{equation}{0}

\section{Coefficients $M_{j}$}\label{M}

To evaluate the coefficients
\beq
 M_j=\int_0^1
 \left(\tan^{2\Gr_1-1}\fr{\pi x}{2}+\cot^{2\Gr_1-1}\fr{\pi x}{2}\right)
 T_j(\cos\pi x)dx, \quad j=0,1,\ldots,
 \label{C.1}
 \eeq
we make the substitution $\Gz=\cos\pi x$ and  rewrite this expression as
\beq
M_j=\fr{1+(-1)^j}{\pi}\int_{-1}^1 (1-\Gz)^{\Gr_1-1}(1+\Gz)^{-\Gr_1}T_j(\Gz)d\Gz.
\label{3.11}
\eeq
It follows immediately that $M_1=M_3=\ldots=0$. To compute
the integral in the even case, $j=2m$, we
consider the integral
\beq
\int_{-1}^1(1-\Gz)^{\Ga_1}(1+\Gz)^{\Ga_2} T_j(\Gz)d\Gz,
\label{3.11'}
\eeq
make the substitution $\Gz=2t-1$ and
express the Chebyshev polynomials of the first kind through the 
hypergeometric function
\beq
T_j(2t-1)=F(-j,j; 1/2; 1-t).
\label{3.12}
\eeq
Then we change the order of integration and summation,  evaluate the new integrals in terms of the $\GG$-functions and have
$$
\int_{-1}^1(1-\Gz)^{\Ga_1}(1+\Gz)^{\Ga_2} T_j(\Gz)d\Gz=
\fr{2^{\Ga_1+\Ga_2+1}\GG(\Ga_1+1)\GG(\Ga_2+1)}{\GG(\Ga_1+\Ga_2+2)}
$$
\beq
\times {}_3 F_2\left(
-j,  j, \Ga_1+1;
1/2, \Ga_1+\Ga_2+2;  1\right), \quad \R\Ga_1>-1, \quad \R\Ga_2>-1.
\label{3.13}
\eeq
Alternatively, this result   can be derived from formula
16.4 (3) in {\bf(\ref{bat2})}.

Note that the formulas 16.1 (2) in {\bf(\ref{bat2})} and  7.347 (1) in {\bf(\ref{gra})}
corresponding to (\ref{3.13})
have the same error. In addition,   formulas 16.1 (1), 16.1 (21) and  16.1 (22) in {\bf(\ref{bat2})} and
the corresponding to  16.1 (22) formula  7.347 (2) in {\bf(\ref{gra})} need  
also to be corrected. They should read, respectively,
$$
\int_{-1}^1(1-\Gz)^{-1/2}(1+\Gz)^{\Ga} T_j(\Gz)d\Gz=
\fr{2^{\Ga+1/2}\sqrt{\pi}\GG(\Ga+1)\GG(\Ga+3/2)}{\GG(\Ga+3/2+j)\GG(\Ga+3/2-j)}.
$$
$$
\int_{-1}^1(1-\Gz)^{1/2}(1+\Gz)^{\Ga} U_j(\Gz)d\Gz=
\fr{2^{\Ga+1/2}\sqrt{\pi}\GG(\Ga+1)\GG(\Ga+1/2)(j+1)}{\GG(\Ga+5/2+j)\GG(\Ga+1/2-j)},
$$
$$
\int_{-1}^1(1-\Gz)^{\Ga_1}(1+\Gz)^{\Ga_2} U_j(\Gz)d\Gz=
\fr{2^{\Ga_1+\Ga_2+1}\GG(\Ga_1+1)\GG(\Ga_2+1)(j+1)}{\GG(\Ga_1+\Ga_2+2)}
$$
\beq
\times {}_3 F_2\left(
-j,  j+2, \Ga_1+1;
3/2, \Ga_1+\Ga_2+2;  1\right), \quad \R\Ga_1>-1, \quad \R\Ga_2>-1.
\label{3.14}
\eeq
Here, $U_j(\Gz)$ are  the Chebyshev polynomials of the second kind. On employing now formula (\ref{3.13}) we obtain that $M_{2m}$ is a finite sum given by
\beq
 M_{2m}=\fr{2}{\sin\pi\Gr_1}\sum_{j=0}^{2m}\fr{(-2m)_j(2m)_j(\Gr_1)_j}{(1/2)_j(j!)^2},
 \quad m=0,1,\ldots,
 \label{C.2}
\eeq
and, in particular, $M_0=2\csc\pi\Gr_1$.

\setcounter{equation}{0}

\section{Complete singular integral equation with the Cauchy kernel}\label{Cauchy}

We aim to show that it is not surprising that the solvability condition of
equation (\ref{3.1}) coincides with the first equation (\ref{3.33}) 
of the infinite
system of algebraic equations (\ref{3.31}). The same result can be derived for the classical equation with the Cauchy kernel
\beq
\fr{1}{\pi}\int_0^1\left[\fr{1}{\Gx-x}+K(x,\Gx)\right]\Gf(\Gx)d\Gx=C-F(x),\quad 0<x<1. 
\label{D.1}
\eeq
in the class of H\"older functions bounded at the endpoints. It is known {\bf(\ref{gak})} that this equation is solvable if and only if the constant $C$
is chosen to be
\beq
C=\fr{1}{\pi}\int_0^1\fr{F(x)dx}{\sqrt{x(1-x)}}+\fr{1}{\pi^2}\int_0^1\int_0^1
\fr{K(x,\Gx)\Gf(\Gx)d\Gx dx}{\sqrt{x(1-x)}}.
\label{D.2}
\eeq
Then the characteristic equation ($K(x,\Gx)\equiv 0$), $\hat S[\Gf](x)=C-F(x)$, $0<x<1$,
 in the class of functions chosen, admits the unique solution $\Gf(x)=-\hat\CS^{-1}[F](x)$, where the inverse operator is given by
 \beq
 \hat\CS^{-1}[F](x)=-\fr{\sqrt{x(1-x)}}{\pi}\int_0^1\fr{F(\Gx)d\Gx}{\sqrt{\Gx(1-\Gx)}(\Gx-x)}, \quad 0<x<1,
 \label{D.2'}
 \eeq
 and $\hat\CS^{-1}[1]=0$, $0<x<1$ (recall that we derived the same result (\ref{3.27''}) for the operator $\CS^{-1}$).

 Now, on returning to the complete equation with the Cauchy kernel, 
 we expand the unknown function $\Gf(x)$ as
\beq
\Gf(x)=\sqrt{x(1-x)}\sum_{j=0}^\infty b_jU_j(2x-1)
\label{D.3}
\eeq
and employ the spectral relation for the Chebyshev polynomials of the second kind
\beq
\int_0^1\fr{\sqrt{\Gx(1-\Gx)}}{\Gx-x}U_j(2\Gx-1)d\Gx=-\fr{\pi}{2}T_{j+1}(2x-1), \quad 0<x<1, \quad j=0,1,\ldots.
\label{D.4}
\eeq
Then because of the orthogonality of the Chebyshev polynomials, we reduce the integral
equation (\ref{D.1}) to
the following system of linear algebraic equations:
\beq
\sum_{j=0}^\infty k_{0j}b_j=C-f_0
\label{D.5}
\eeq
and
\beq
-\fr{b_n}{4}+\sum_{j=0}^\infty k_{nj}b_j=-f_n, \quad n=1,2,\ldots.
\label{D.6}
\eeq
Here,
$$
 k_{nj}=\fr{1}{\pi^2}
 \int_0^1
 \int_0^1 
 K(x,\Gx)\sqrt{\Gx(1-\Gx)}U_j(2\Gx-1)
 \fr{T_n(2x-1)d\Gx dx}{\sqrt{x(1-x)}},
$$
\beq
f_n=\fr{1}{\pi}\int_0^1 \fr{F(x)T_{n}(2x-1)dx}{\sqrt{x(1-x)}}.
\label{D.7}
\eeq
It becomes evident, upon substituting the series (\ref{D.3}) into the solvability condition 
(\ref{D.2}) and using the notations (\ref{D.7}), that equations (\ref{D.2})
and (\ref{D.5}) are equivalent.

\end{document}